\documentclass[reqno]{amsart}
\usepackage{amsmath}%
\usepackage{amsfonts}%
\usepackage{amssymb}%
\usepackage{graphicx}
\usepackage{mathrsfs}
\usepackage{hyperref}
\usepackage{fullpage}
\usepackage{lineno}
\usepackage{setspace}
\usepackage{enumitem}
\usepackage{xcolor}
\usepackage{pgfplots}
\usepackage{comment}

\newtheorem{theorem}{Theorem}
\theoremstyle{plain}

\newtheorem{claim}[theorem]{Claim}
\newtheorem{construction}[theorem]{Construction}

\newtheorem{corollary}[theorem]{Corollary}

\newtheorem{fact}[theorem]{Fact}
\newtheorem{lemma}[theorem]{Lemma}

\newtheorem{proposition}[theorem]{Proposition}

\numberwithin{equation}{section}
\numberwithin{theorem}{section}

\pgfdeclarelayer{bg}    
\pgfsetlayers{bg,main}  

\def\ba{\backslash}
\def\A{\mathcal{A}}
\def\B{\mathcal{B}}
\def\C{\mathcal{C}}

\def\F{\mathcal{F}}
\def\G{\mathcal{G}}
\def\h{\mathcal{H}}
\def\K{\mathcal{K}}
\def \L{\mathcal{L}}
\def\M{\mathcal{M}}

\def\Q{\mathcal{Q}}
\def\R{\mathcal{R}}
\def\T{\mathcal{T}}
\def\S{\mathcal{S}}

\def\pp{\mathcal{P}}
\def \a{\alpha}
\def\epsilon{\varepsilon}
\def \e{\epsilon}
\def \ext{\mathrm{ext}}
\def \even{\text{even}}
\def \odd{\text{odd}}
\def \r{\gamma}

\begin{document}


\title[On Hamilton $(k/2)$-cycles]{Minimum degree thresholds for Hamilton $(k/2)$-cycles in $k$-uniform hypergraphs}
\date{\today}
\author{Hi\d{\^{e}}p H\`an}
\author{Jie Han}
\author{Yi Zhao}

\address{Departamento de  Matem\'atica y Ciencia de la Computaci\'on, Universidad de Santiago de Chile, Las Sophoras 173, Santiago, Chile}
\thanks{HH is supported by the FONDECYT Regular grant 1191838.}
\email[Hiep Han]{hiep.han@usach.cl}
\address
{Department of Mathematics, University of Rhode Island, 5 Lippitt Road, Kingston, RI, 02881
}
\email[Jie Han]{jie\_han@uri.edu}
\address
{Department of Mathematics and Statistics, Georgia State University, Atlanta, GA 30303}
\thanks{JH is partially supported by Simons Collaboration Grant for Mathematicians \#630884. Part of this research was done while JH was a PhD student at Georgia State University under the supervision of YZ}
\thanks{YZ is partially supported by NSF grants DMS-1400073 and DMS-1700622.}
\email[Yi Zhao]{yzhao6@gsu.edu}
\date{\today}
\begin{abstract}
For any even integer $k\ge 6$, integer $d$ such that $k/2\le d\le k-1$, and sufficiently large $n\in (k/2)\mathbb N$, we find a tight minimum $d$-degree condition that guarantees the existence of a Hamilton $(k/2)$-cycle in every $k$-uniform hypergraph on $n$ vertices.
 When $n\in k\mathbb N$, the degree condition coincides with the one for the existence of perfect matchings provided by R\"odl, Ruci\'nski and Szemer\'edi (for $d=k-1$) and Treglown and Zhao 
 (for $d\ge k/2$), and thus our result strengthens theirs in this case.
\end{abstract}

\maketitle


\section{Introduction}\label{sec:intro}
The study of Hamilton cycles is an important topic in graph theory with a long history. 
In recent years, researchers have worked on extending the classical theorem of Dirac on Hamilton cycles to hypergraphs 
and we refer to~\cite{BHS, GPW, HZ1, RR, BMSSS1, BMSSS2, RRRSS, HZ_forbidHC} for some recent results and to \cite{KuOs14ICM, RR, zsurvey} for surveys  on this topic.
Given $k\ge 2$, a $k$-uniform hypergraph (in short, \emph{$k$-graph}) consists of a vertex set $V$ and an edge set $E\subseteq \binom{V}{k}$, where $\binom{V}{k}$
denotes the family of all  $k$-element subsets of~$V$. 
Given a $k$-graph $\h$ with a set $S$ of $d$ vertices, where $1 \le d \le k-1$, we define $\deg_{\h} (S)$ to be the number of edges containing $S$ (the subscript $\h$ is omitted if it is clear from the context). The \emph{minimum $d$-degree $\delta _{d} (\h)$} of $\h$ is the minimum of $\deg_{\h} (S)$ over all $d$-vertex sets $S$ in $\h$.  We refer to  $\delta _1 (\h)$ as the \emph{minimum vertex degree} and  $\delta _{k-1} (\h)$ the \emph{minimum codegree} of $\h$.
For $1\le \ell< k$,
a $k$-graph is called an \emph{$\ell$-cycle} if its vertices can be ordered cyclically such that each of its edges consists of $k$ consecutive vertices and every two consecutive edges (in the natural order of the edges) share exactly $\ell$ vertices. In $k$-graphs, a $(k-1)$-cycle is often called a \emph{tight} cycle.
We say that a $k$-graph contains a \emph{Hamilton $\ell$-cycle} if it contains an $\ell$-cycle as a spanning subhypergraph. Note that a Hamilton $\ell$-cycle of a $k$-graph on $n$ vertices contains exactly $n/(k - \ell)$ edges, implying that $k- \ell$ divides $n$.

Confirming a conjecture of Katona and Kierstead \cite{KK}, R\"odl, Ruci\'nski and Szemer\'edi \cite{RRS06, RRS08} showed that for any fixed $k$, every $k$-graph $\h$ on $n$ vertices with $\delta_{k-1}(\h)\ge n/2 + o(n)$ contains a tight Hamilton cycle. 
When $k-\ell$ divides both $k$ and $|V|$, a tight cycle on $V$ contains an $\ell$-cycle on $V$. 
Thus the result in \cite{RRS08} implies that for all $1\le \ell<k$ such that $k- \ell$ divides $k$, every $k$-graph $\h$ on $n\in (k-\ell)\mathbb{N}$ vertices with $\delta_{k-1}(\h)\ge n/2 + o(n)$ contains a Hamilton $\ell$-cycle. 
This is best possible up to the $o(n)$ term by a construction given by Markstr\"om and Ruci\'nski~\cite{MaRu}. 
R\"odl, Ruci\'nski and Szemer\'edi~\cite{RRS11} eventually determined the minimum codegree threshold for tight Hamilton cycles in 3-graphs for sufficiently large $n$, which is $\lfloor n/2\rfloor$.  

After a series of efforts~\cite{KKMO, HS, KMO}, the minimum codegree conditions for $\ell$-Hamiltonicity were determined asymptotically.
R\"odl and Ruci\'nski \cite[Problem 2.9]{RR} raised the question concerning  the \emph{exact} minimum codegree condition for $\ell$-Hamiltonicity when $n$ is sufficiently large. In the rest of the paper, unless stated otherwise, we assume that~$n$ is sufficiently large.
The case $k=3$ and $\ell=1$ was solved by Czygrinow and Molla \cite{CzMo}, and the last two authors~\cite{HZ2} determined this threshold for all $k\ge 3$ and $\ell<k/2$. 
Recently, the case $k=4$ and $\ell=2$ was determined by Garbe and Mycroft~\cite{GaMy}.
We continue this line of research and obtain the minimum codegree threshold for $\ell$-Hamiltonicity with even $k\ge 6$ and $\ell=k/2$.

\begin{theorem}\label{thm_codegree}
For all even integers $k\ge 6$ there exists $n_0$ such that  the following holds for every $n\in \frac k2 \mathbb{N}$ with $n\ge n_0$.
If $\h$ is a $k$-uniform hypergraph on $n$ vertices such that 
\[
\delta_{k-1}(\h)> \begin{cases}
 n/2 -k+1 &\text{ if $n\in k\mathbb N$ and $n/2 - n/k$ is even}\\
\lfloor n/2\rfloor -k+2 &\text{ otherwise, }\\
\end{cases}
\]
then $\h$ contains a Hamilton $(k/2)$-cycle.
\end{theorem}

We note that by the result of~\cite{GaMy}, Theorem~\ref{thm_codegree} also holds for $k=4$ with the same minimum codegree threshold.
Theorem~\ref{thm_codegree} follows from our main result, Theorem~\ref{thmmain}, which goes far beyond the minimum codegree condition and determines the minimum $d$-degree condition for $(k/2)$-Hamiltonicity for every $k/2\le d\le k-1$.
To state our main result, we first  introduce some notation.

\subsection{Lower bound constructions and main result}
Our constructions will build upon the ones for perfect matchings from~\cite{TrZh12}.
In fact, our extremal hypergraphs coincide with the ones in~\cite{TrZh12} when $n\in k\mathbb{N}$. 
Let a vertex set $V$  with a partition $V= A\dot\cup B$ be given. 
A set $S\subseteq V$ is \emph{odd} (\emph{w.r.t.}~$A$) or \emph{even} (\emph{w.r.t.}~$A$) if it intersects $A$ in an odd or even number of vertices, respectively.
Let $E_{\odd}(A,B)$ denote the family of all odd $k$-element subsets of $V$ and let  $E_{\even}(A,B)$ be the family of even $k$-element subsets of $V$.
Define $\B_{n,k}(A,B)$ and $\overline \B_{n,k}(A,B)$, respectively, to be the $k$-graph with vertex set $V=A\dot\cup B$ and edge set 
$E_{\odd}(A,B)$ and $E_{\even}(A,B)$, respectively. 
A star $\mathcal S_{n,k}$ is an $n$-vertex $k$-graph that consists of all $k$-sets containing a fixed vertex $v$. Let $\B'_{n,k}(A, B)$ be the hypergraph obtained from $\B_{n,k}(A, B)$ by adding a star $\mathcal S_{|A|, k}$ into $A$.

\begin{construction}
\label{cons}
Given an even integer $k\ge 4$ and an integer $n\in  \frac k2\mathbb N$, our family  $\h_{\ext}(n,k)$ of extremal $k$-graphs is defined as follows. 
\begin{itemize}
\item For $n\in k\mathbb N$ the family $\h_{\ext}(n,k)$ contains all hypergraphs $\B_{n,k}(A,B)$ where $n/k-|A|$ is odd and all  $\overline{\mathcal B}_{n,k}(A,B)$ where $|A|$ is odd. 
\item For $n\in \frac k2\mathbb N\setminus k\mathbb N$,  $\h_{\ext}(n,k)$ contains all hypergraphs $\B_{n,k}(A, B)$ together with 
\begin{itemize}
\item all $\B'_{n,k}(A,B)$ when $k\in 4\mathbb N$ and $\lfloor n/k \rfloor - |A|$ is odd; and
\item all $\B'_{n,k}(A,B)$ when $k\in 2\mathbb N\setminus 4\mathbb N$ and $\lfloor n/k \rfloor - |A|$ is even.
\end{itemize}
\end{itemize}
\end{construction}
For $n\in k\mathbb N$, it was shown in~\cite{TrZh12} that any hypergraph in $\h_{\ext}(n,k)$ contains no perfect matching and thus no Hamilton $(k/2)$-cycle (because a Hamilton $(k/2)$-cycle consists of two disjoint perfect matchings). We will show that no hypergraph in $\h_{\ext}(n,k)$ contains a Hamilton $(k/2)$-cycle when $n\in \frac k2\mathbb N\setminus k\mathbb N$.
To do so we will represent a Hamilton $(k/2)$-cycle $\C$ as a sequence of pairwise disjoint $(k/2)$-sets 
$L_1, \dots, L_t$ with $t= 2n/k$ such that $L_i\cup L_{i+1}\in E(\C)$ for all $i\in [t]$, where $L_{t+1}:=L_1$. 
Further,  we associate to $\C$ the binary string $b_1b_2\cdots b_t$, called the \emph{binary representation of $\C$} (\emph{w.r.t.}~$A$), 
defined by $b_i=0$ if $|L_i\cap A|$ is even and  $b_i=1$ otherwise.

\begin{proposition}\label{prop:lowerbounds}
No $k$-graph in $\h_{\ext}(n,k)$ contains a Hamilton $(k/2)$-cycle. 
\end{proposition}
\begin{proof}
For $n\in k\mathbb N$ note that no $k$-graph in $\h_{\ext}(n,k)$ contains a perfect matching (and thus none contains a Hamilton $(k/2)$-cycle either).
This is because in $\overline{\mathcal B}_{n,k}(A,B)$ all edges are even while $|A|$ is odd,
and in  ${\mathcal B}_{n,k}(A,B)$ all edges are odd while the cardinality of a perfect matching,  $n/k$,  and $|A|$ have different parities.

Consider now $n\in \frac k2\mathbb N\setminus k\mathbb N$ and suppose that some $k$-graph in $\h_{\ext}(n,k)$
contains a Hamilton $(k/2)$-cycle~$\C$. Let 
 $b_1b_2\dots b_t$ be  the binary representation of $\C$  and note that an odd edge in $\C$ corresponds to a 01 or 10 in this representation.
Thus $\C$ cannot consists of odd edges only, as 
then  $b_i\neq b_{i+1}$ holds for all $i\in [t]$ ($b_{t+1}=b_1$) yet $t=2{n}/{k}$ is odd. This implies that $\B_{n,k}(A, B)$ contains no Hamilton $(k/2)$-cycles.
Moreover, any such cycle $\C$ in $\B'_{n,k}(A, B)$ must contain at least one even edge. 

We claim that every Hamilton $(k/2)$-cycle $\C$ in $\B'_{n,k}(A, B)$ contains exactly one even edge.
To this end note that the even edges  in $\B'_{n,k}(A, B)$ form a star, thus  $\C$ cannot contain three or more even edges, as there would be
two disjoint ones otherwise. Furthermore, as the star is entirely contained in $A$, an even edge in $\C$
corresponds to a pair 00 when $k\in 4\mathbb N$ (and thus $k/2$ is even) and to a pair 11 when $k\in 2\mathbb N\setminus 4\mathbb N$  (and thus $k/2$ is odd).
In the first case we conclude that the number of odd edges is twice the number of 1-entries in the binary representation, because 
each odd edge gives rise to exactly one 1-entry while one such entry is a witness for two odd edges.
Similarly, the number of odd edges is twice the number of 0-entries
in the second case. Therefore the number of odd edges in $\C$ is even while $t=2n/k$ is odd, which implies that the number of even edges in $\C$ is one, as claimed.

We conclude that for $k\in 4\mathbb N$ the cycle $\C$ has the form $00101\cdots01$, thus
contains $\lfloor n/k \rfloor$ odd $(k/2)$-sets. However, $\lfloor n/k \rfloor$ and $|A|$ have different parities 
which yields a contradiction.
For $k\in 2\mathbb N\setminus 4\mathbb N$ the cycle $\C$ has the form $11010 \cdots 10$ and
contains therefore $\lfloor n/k \rfloor+1$ odd $(k/2)$-sets. This implies that $|A|$ and  $\lfloor n/k \rfloor+1$ have the same parity
which  yields a contradiction to the assumption that $\lfloor n/k \rfloor - |A|$ is even. 
\end{proof}

The following is our main result, which states that $k$-graphs with minimum $d$-degree larger than the ones in $\h_{\ext}(n,k)$ must contain a Hamilton $(k/2)$-cycle.
Given positive integers $d<k\le n$ such that $k$ is even and $k/2$ divides $n$, 
let $\overline\delta(n,k,d)$ be the maximum of the minimum $d$-degree among all the hypergraphs in $\h_{\ext}(n,k)$. 

\begin{theorem}[Main Result] \label{thmmain}
For even integers $k\ge 6$, $k/2\le d\le k-1$ and sufficiently large integer $n\in \frac k2\mathbb{N}$ the following holds. 
Suppose $\h$ is a $k$-graph on $n$ vertices satisfying $\delta_{d}(\h) > \overline\delta(n,k,d)$,
then $\h$ contains a Hamilton $(k/2)$-cycle.
\end{theorem}

When $k$ is even, it is easy to see that $\delta_{k-1}(\overline\B_{n,k}(A, B))=\min\{|A|-k+1, |B|-k+1\}$ and $\delta_{k-1}(\B_{n,k}(A, B))=\delta_{k-1}(\B'_{n,k}(A, B))=\min\{|A|-k+2, |B|-k+2\}$.
Thus, it is straightforward to check that 
\[
\overline\delta(n,k,k-1)= \begin{cases}
 n/2 -k+1 &\text{ if $n\in k\mathbb N$ and $n/2 - n/k$ is even}\\
\lfloor n/2\rfloor -k+2 &\text{ otherwise. }\\
\end{cases}
\]
Theorem~\ref{thm_codegree} is therefore a special case of Theorem~\ref{thmmain}. 
Moreover, given positive integers $d<k\le n$ such that $k$ divides $n$ ($k$ is not necessarily even),
let $\delta(n,k,d)$  be the maximum of the minimum $d$-degree among all the hypergraphs from the first class of Construction~\ref{cons}.
Then $\overline\delta(n,k,d)= \delta(n,k,d)$ when $k$ is even and $n\in k\mathbb N$.
Extending a result of R\"odl, Ruci\'nski and Szemer\'edi~\cite{RRS09}, 
Treglown and Zhao \cite{TrZh13} showed that if $\delta_{d}(\h) > \delta(n,k,d)$, then every $n$-vertex $k$-graph $\h$ contains a perfect matching.
Theorem~\ref{thmmain} shows that, for even $k\ge 6$, the minimum $d$-degree that forces the existence of a perfect matching actually forces a Hamilton $(k/2)$-cycle, a union of two disjoint perfect matchings. Therefore Theorem~\ref{thmmain} strengthens the results of \cite{RRS09,TrZh13}.

We note that, however, the precise values of $\delta(n,k,d)$ and $\overline\delta(n,k,d)$ when $d\le k-2$ are only known to be $(1/2+o(1))\binom{n-d}{k-d}$, 
see  \cite{TrZh12} for details.

\subsection{Proof of Theorem \ref{thmmain}}

As a common approach to obtain exact results, Theorem~\ref{thmmain} is proven by distinguishing an \emph{extremal} case from a \emph{nonextremal} case and solve them separately.
Let $\e>0$ and suppose that $\h$ and $\h'$ are $k$-graphs on $n$ vertices. 
We say that $\h$ is \emph{$\e$-close to $\h'$}, and write $\h=\h'\pm \e n^k$, if $\h$ can be made a copy of $\h'$ by adding and deleting at most $\e n^k$ edges.
Suppose that  $\h$ is a $k$-graph  with minimum $d$-degree 
$\delta_{d}(\h) \ge (\frac 12-o(1)) \binom{n-d}{k-d}$ and $o(1)$-close to some $k$-graph in $\h_{\ext}(n,k)$, then~$\h$ must be $o(1)$-close to some 
$\B_{n,k}(A,B)$ or $\overline \B_{n,k}(A,B)$ with $|A|= |B|= n/2$ as well. 
In the following we simply write~$\B_{n,k}$ and~$\overline \B_{n,k}$ to indicate that there is an implicit partition $A\cup B$ of equal size.

\begin{theorem}[Nonextremal Case]\label {lemNE}
For any integer $k\ge 4$ even, $k/2\le d\le k-1$ and $\e>0$ there exist $\r>0$ and $n_{\ref{lemNE}}$ such that for every $k$-graph $\h=(V, E)$ on $n\ge n_{\ref{lemNE}}$ vertices with $n\in (k/2)\mathbb N$ the following holds. Suppose that $\h$ is not $\e$-close to any $\B_{n,k}$ or $\overline \B_{n,k}$
and $\delta_{d}(\h) \ge (\frac 12-\r) \binom{n-d}{k-d}$, then $\h$ contains a Hamilton $(k/2)$-cycle.
\end{theorem}

\begin{theorem}[Extremal Case]\label {lemE}
For any integer $k\ge 6$ even and $k/2\le d\le k-1$, there exist $\e>0$ and $n_{\ref{lemE}}\in \mathbb N$ such that for every $k$-graph $\h=(V, E)$ on $ n\ge n_{\ref{lemE}}$ vertices with $n\in (k/2)\mathbb N$ the following holds. Suppose that $\delta_d(\h)>\overline\delta(n,k,d)$ and $\h$ is $\e$-close to a $\B_{n,k}$ or a $\overline \B_{n,k}$, then $\h$ contains a Hamilton $(k/2)$-cycle.
\end{theorem}

Theorem \ref{thmmain} follows from Theorems \ref{lemE} and \ref{lemNE} immediately by choosing $\e$ from Theorem~\ref{lemE} and letting $n_{\ref{thmmain}}=\max\{n_{\ref{lemE}}, n_{\ref{lemNE}}\}$.

Let us briefly discuss our proof ideas.
Theoreom~\ref{lemNE} is proven in Section~\ref{sec:3}.  
Following previous work \cite{RRS06, RRS08, RRS11, HS, KMO, BHS}, we use the absorbing method initiated by R\"odl, Ruci\'nski and Szemer\'edi. 
More precisely, we find the desired Hamilton cycle by three lemmas: the Absorbing Lemma  (Lemma~\ref{lemA}), the Reservoir Lemma (Lemma~\ref{lemR}), and the Path-cover Lemma (Lemma~\ref{lemP}).
In fact, both the Reservoir Lemma and Absorbing Lemma can be easily derived from a Connecting Lemma (Lemma~\ref{lemC}), which says that either $\h$ is extremal or any two $(k/2)$-sets in $\h$ must have many sets that `connect' them as a $(k/2)$-path. 
To prove the Path-cover Lemma, we slightly strengthen a result of Markstr\"om and Ruci\'nski~\cite{MaRu} on matchings in $k$-graphs and use the regularity method to obtain an almost path-cover of $\h$.
The main technicality lies in the proof of the Connecting Lemma, in which we follow the stability method along a scheme given by Treglown and Zhao~\cite{TrZh12, TrZh13}.
The proof of Theorem~\ref{lemE} is more challenging with one of 
the main complications stemming from the fact that there are several extremal $k$-graphs for the 
problem and different strategies must be used to overcome the (parity) constraint in each case (see Section~\ref{sec:8} for a more detailed outline).
Suppose $\h$ satisfies $\delta_d(\h)>\overline\delta(n,k,d)$ and is close to $\B_{n,k}$ or $\overline \B_{n,k}$.
Using the minimum degree condition,  we can build a short path which can break the parity barriers and be extended to a Hamilton cycle of $\h$. 
The argument of constructing this short path crucially relies on Lemmas~\ref{lem:deg_inter},~\ref{prop:e123} and \ref{thm:deg_bridge} from below, three results  concerning $k$-graphs with forbidden intersections.
These lemmas belong to a  line of research which is central in extremal set theory with a long and influential history. We feel that they are of independent interest and may find applications beyond the one considered here. Therefore we will discuss the two lemmas in more detail  in the following subsection.

\subsection{Breaking the parity barriers and set systems with forbidden intersections}

To break the parity barriers one is of course interested in the existence of additional even
edges in case of $\B_{n,k}(A, B)$ and $\B_{n,k}'(A, B)$  and odd edges in case of $\overline\B_{n,k}(A, B)$ respectively.  
Indeed, in the extremal case of \cite{TrZh12}, the existence of one such edge is enough to overcome the extremal examples 
$\B_{n,k}(A, B)$ and $\overline\B_{n,k}(A, B)$ 
and find the perfect matching. 
As shown by the following example, our problem is more complicated, as additional edges may not be enough to provide a Hamilton $(k/2)$-cycle.
Let us first discuss the case $\overline\B_{n,k}(A, B)$ for 
$n\in k\mathbb{N}$ and $V(\h)=A\cup B$ with odd $|A|$. Assume that $\h$ consists of $\overline\B_{n,k}(A, B)$ together with a set of odd edges such that
\begin{equation}\label{eq:0k/2}
\text{no  two odd edges are disjoint or intersect in exactly $k/2$ vertices.}
\end{equation}
By Proposition~\ref{prop:lowerbounds} we know that $\overline\B_{n,k}(A, B)$ contains no Hamilton $(k/2)$-cycle so a possible Hamilton $(k/2)$-cycle~$\C$ in~$\h$ would have to use odd edges.
But due to~\eqref{eq:0k/2}, $\C$ must contain exactly one odd edge, which is impossible by considering the binary representation of $\C$.

Concerning  the parity barriers posed by $\B_{n,k}(A, B)$ and $\B'_{n,k}(A, B)$ recall that the latter itself consists of $\B_{n,k}(A, B)$ with a star $\S_{|A|,k}$ added to $A$.
The proof of Proposition~\ref{prop:lowerbounds} essentially showed that a $k$-graph~$\h$ (e.g., $\B'_{n,k}(A, B)$) still contains no Hamilton $(k/2)$-cycle, if it 
consists of $\B_{n,k}(A, B)$ together with a set of even edges in $A$ such that 
\begin{equation}\label{eq:no3even}
\text{no three of them are part of a $(k/2)$-path.}
\end{equation}  

The following result will be crucial for overcoming the barriers mentioned above.
It bounds the number of edges in a $k$-graph with a specific forbidden intersection pattern and we derive it from a result of Frankl and F\"uredi~\cite{FrFu85}.

\begin{lemma}\label{lem:deg_inter}
For every  even $k\ge 4$  there is a $c>0$ such that the following holds. 
Suppose $\h$ is a $k$-graph on $n$ vertices such that $|e_1\cap e_2|\neq 0, k/2$ for any two edges $e_1, e_2$ in $\h$. 
Then $e(\h)<c n^{k/2-1}$.
\end{lemma}

\begin{proof}
We first recall a classical theorem by Frankl and F\"uredi~\cite{FrFu85} concerning ``forbidding just one intersection''. It states that for any $1\le \ell\le k'-1$ there is a $c>0$ 
so that the following holds. 
\begin{equation}\label{eq:FF}
\text{If  $\mathcal F\subset\tbinom{[n]}{k'}$ is  such that $|A\cap B|\neq \ell$ for all $A,B\in \mathcal F$, then  $|\mathcal F|\leq c n^{\max\{\ell,k'-\ell-1\}}$.}
\end{equation}

Turning to the proof of the lemma consider an arbitrary edge $\{v_1, \dots, v_k\}\in E(\h)$ and let $N(v)=\{e\setminus\{v\}\colon e\in E(\h)\}$. 
Since $\h$ is intersecting (i.e., any two of its members have a non-empty intersection), 
any edge of $\h$ intersects $\{v_1, \dots, v_k\}$ and thus
$e(\h)\leq \sum_{i\in[k]} |N(v_i)|$.
 By the assumption on $\h$ 
we know that $N(v_i)$ is  $(k-1)$-uniform and the intersection of any two edges of $N(v_i)$
has size distinct from $k/2-1$. Thus applying~\eqref{eq:FF}
with $k'=k-1$ and $\ell=k/2-1$ on each $N(v_i)$ we obtain $e(\h)\le \sum_{i\in [k]}|N(v_i)|\le kc_{k-1}' n^{k/2-1}$.
\end{proof}

By applying Lemma~\ref{lem:deg_inter} and the Hilton-Milner theorem~\cite{HM67}, we obtain the following lemma and will use it to address the barrier \eqref{eq:no3even}. 

\begin{lemma}\label{prop:e123}
Let $k\ge 4$ be even and $n$ be sufficiently large. 
Suppose $\h$ is a $k$-graph on $n$ vertices such that $e(\h) \geq 2k^2\binom{n-2}{k-2}$ and $\h$ is not a subgraph of $\mathcal S_{n,k}$. 
Then there exist three edges $e_1, e_2, e_3$ such that $e_1\cap (e_2\cup e_3)=\emptyset$ and $|e_2\cap e_3|\in \{0, k/2\}$.
In particular, the conclusion holds for $k$-graphs $\h$ such that $\delta_d(\h) > \delta_d(\mathcal S_{n,k})$ for any $d<k$.
\end{lemma}

\begin{proof}
A classical result of Hilton and Milner~\cite{HM67} states that if $\h$ is intersecting but not a subgraph of~$\mathcal S_{n,k}$, then $|E(\h)|\le k\binom{n-2}{k-2}$.
This together with our assumptions implies 
that~$\h$ is not intersecting.
Let  $e_1$ and $e_2$ be  two disjoint edges of $\h$ and 
let~$\h'$ denote the subgraph obtained from~$\h$ by removing all edges intersecting both $e_1$ and $e_2$. 
Then $e(\h')\geq e(\h)-k^2\binom{n-2}{k-2}> k^2\binom{n-2}{k-2}$ and  
we may assume that $\h'$ contains no edge which is disjoint from $e_1\cup e_2$ since we would be done otherwise. 
Then $\h'$ can be partitioned into $\h_1$ and $\h_2$, where $\h_i$ contains all edges intersecting $e_i$ for $i=1,2$ (thus not intersecting $e_{3-i}$). 
We have $e(\h_i)> \frac{k^2}2\binom{n-2}{k-2}$ for some $i\in[2]$ and applying Lemma \ref{lem:deg_inter} (note that $k-2> \frac k2 -1$) we obtain the desired third edge $e_3$.

Note that if $\delta_d(\h)> \delta_d(\mathcal S_{n, k})=  \binom{ n - 1 - d}{k-1 - d}$, then $\h$ is not a subgraph of $\mathcal S_{n,k}$, and 
\[
e(\h) > \binom{n}{d} \binom{ n - 1 - d}{k-1 - d} / \binom kd = \Omega(n^{k-1}) \geq 2k^2\binom{n-2}{k-2}.
\]
So the second part of the lemma follows.
\end{proof}

By Lemma~\ref{lem:deg_inter}, we also obtain the following lemma and will use it to address the barrier \eqref{eq:0k/2}. 
\begin{lemma}\label{thm:deg_bridge}
Given an even integer $k\ge 4$ and an integer $d\ge k/2$, let $n$ be sufficiently large.
Suppose $\h$ is an $n$-vertex $k$-graph with a partition $V(\h)=A\cup B$ such that $|A|, |B|\ge 0.4 n$ and $\delta_d(\h) > \delta_d(\B_{n, k}(A, B))$ (respectively, $\delta_d(\h) > \delta_d(\overline \B_{n, k}(A, B))$). 
Then $ \h\cap \overline{\B}_{n, k}(A, B)$ (respectively, $\h\cap \B_{n, k}(A, B)$) contains two edges $e_1, e_2$ such that $|e_1\cap e_2|\in \{0, k/2\}$.
\end{lemma}

\begin{proof}
We assume that $\delta_d(\h) > \delta_d(\B_{n, k}(A, B))$ because the other case can be proved similarly.
Let $E_0=E(\h)\cap \overline \B_{n,k}(A, B)$. 
Assume to the contrary that for any two edges $e_1, e_2\in E_0$, we have $|e_1\cap e_2|\neq 0$ or $k/2$. 
Then by Lemma \ref{lem:deg_inter}, $|E_0|=O(n^{k/2-1})$.

On the other hand, we bound $|E_0|$ from below as follows. 
Given the partition $V(\h)=A\cup B$, we partition $\binom{V(\h)}{d}$ into $T_0\cup T_1\cup \cdots\cup T_{d}$ where $T_i=\{X: |X\cap A|=i\}$. 
Note that all $X\in T_i$ have the same degree in ${\B}_{n,k}(A, B)$.
So there exists $j\in \{0,\dots,d\}$ such that the minimum $d$-degree in ${\B}_{n,k}(A, B)$ is achieved by all sets in $T_j$.
Clearly we have
\[
|T_j| =\binom{|A|}{j} \binom{ |B|}{d-j} \ge \binom{0.4 n}{d} = \Omega(n^{k/2})
\]
because $d\ge k/2$.
Moreover, since $\delta_d(\h) > \delta_d(\B_{n, k}(A, B))$ each set in $T_j$ is contained in at least one even edge, thus, we have $|E_0|\ge |T_j| / {\binom kd}= \Omega(n^{k/2})$, a contradiction.
\end{proof}

\subsection*{Notation} Throughout the paper we omit floor and ceiling signs where they do not affect the arguments. Further, we write $\alpha \ll \beta \ll \gamma$ to mean that 
it is possible to choose the positive constants
$\alpha, \beta, \gamma$ from right to left. More
precisely, there are increasing functions $f$ and $g$ such that, given
$\gamma$, whenever we choose some $\beta \leq f(\gamma)$ and $\alpha \leq g(\beta)$, the subsequent statement holds. 
Hierarchies of other lengths are defined similarly.

\section{ Nonextremal Case -- proof of Theorem \ref{lemNE}}\label{sec:3}
In this section we prove Theorem \ref{lemNE}. The following simple and well-known proposition reduces the proof  to the case $d=k/2$. 
\begin{proposition}\label{prop:seq}
Let $0\le d'\le d<k$ and $\h$ be a $k$-graph. If $\delta_d(\h)\ge x\binom{n-d}{k-d}$ for some $0\le x\le 1$, then $\delta_{d'}(\h)\ge x\binom{n-d'}{k-d'}$.\qed
\end{proposition}
The proof of Theorem \ref{lemNE} follows the procedure in \cite{RRS08}. 
A $k$-uniform $\ell$-path of length $t$  is a sequence of vertices $\mathcal P=v_1v_2v_3\cdots v_{(t-1)(k-\ell)+k}$ such that for every $i\in \{0,1,\dots, t-1\}$, $\{v_{i(k-\ell)+1},\dots, v_{i(k-\ell)+k}\}$ forms an edge.
For a $(k/2)$-path $\mathcal P=v_1v_2\cdots v_p$, two $(k/2)$-sets $v_1\cdots v_{k/2}$ and $v_{p-k/2+1}\cdots v_p$ are called the \emph{ends} of $\mathcal P$.
Given a set $C$ and $(k/2)$-sets $A$, $B$, we call $C$ a \emph{connecting $|C|$-set for $A$ and $B$} if $\h[A\cup C\cup B]$ contains a $(k/2)$-path with ends $A$ and $B$.

\begin{lemma}[Reservoir Lemma]\label{lemR}
For an even integer $k\ge 4$ and an integer $d$ with  $k/2\le d\le k-1$, suppose $1/n\ll \r\ll\e\ll 1/k$.
Let $\h=(V,E)$ be a $k$-graph on $n$ vertices with $\delta_{d}(\h) \ge (1/2-\r) \binom{n-d}{k-d}$. 
If $\h$ is not $\e$-close to $\B_{n,k}$ or $\overline \B_{n,k}$, then there is a family $\R$ of $\r^{28} n$ disjoint $(3k/2)$-sets such that 
every pair of $(k/2)$-sets $S, T\subset V(\h)$ has at least $\r^{32}n/3$ connecting $(3k/2)$-sets in $\R$.
\end{lemma}

\begin{lemma}[Absorbing Lemma]\label{lemA}
For an even integer $k\ge 4$ and an integer $d$ with  $k/2\le d\le k-1$, suppose $1/n\ll \r\ll\e\ll 1/k$.
Let $\h=(V,E)$ be a $k$-graph on $n$ vertices with $\delta_{d}(\h) \ge (1/2-\r) \binom{n-d}{k-d}$. 
If $\h$ is not $\e$-close to $\B_{n,k}$ or $\overline \B_{n,k}$, then there exists a $(k/2)$-path $\mathcal P$ in $\h$ with $|V(\mathcal P)|\le 4k\r^{14} n$ such that for all subsets $U\subset V\ba V(\mathcal P)$ of size at most $k\r^{27}n/6$ such that $|U|\in \frac k2\mathbb N$ there exists a $(k/2)$-path $\mathcal Q \subset \h$ with $V(\mathcal Q) = V(\mathcal P)\cup U$ and, moreover, $\mathcal P$ and $\mathcal Q$ have exactly the same ends.
\end{lemma}

\begin{lemma}[Path-cover Lemma]\label{lemP}
For an even integer $k\ge 4$ and an integer $d$ with  $k/2\le d\le k-1$, suppose $1/n\ll 1/p \ll \a \ll \r\ll 1/k$ for some integers $p$ and $n$. 
Let $\h=(V,E)$ be a $k$-graph on $n$ vertices with $\delta_{d}(\h) \ge (1/2-\r) \binom{n-d}{k-d}$. 
Then there is a family of $(k/2)$-paths in $\h$ consisting of at most $p$ paths, which covers all but at most $\a n$ vertices of $\h$.
\end{lemma}

Now we are ready to prove Theorem \ref{lemNE}.
\begin{proof}[Proof of Theorem \ref{lemNE}]
Given an even integer $k\ge 4$ and an integer $d$ with  $k/2\le d\le k-1$, suppose $1/n\ll1/p, \a \ll \r\ll\e\ll 1/k$ for some integers $p$ and $n$.
Let $\h=(V,E)$ be a $k$-graph on $n$ vertices such that $\delta_{d}(\h) \ge (1/2-\r) \binom{n-d}{k-d}$ and assume that $\h$ is not $\e$-close to $\B_{n,k}$ or $\overline \B_{n,k}$.

Since $\h$ is not $\e$-close to $\B_{n,k}$ or $\overline \B_{n,k}$, we can find an absorbing path $\mathcal P_0$ by Lemma \ref{lemA} with ends $S_0, T_0$ and $|V(\mathcal P_0)|\le 4k\r^{14} n$. 
Let $V_1=(V\setminus V(\mathcal P_0))\cup(S_0\cup T_0)$, we claim that $\h[V_1]$ is not $(\e/2)$-close to $\B_{|V_1|,k}$ or $\overline \B_{|V_1|,k}$. 
Suppose instead, that there is a partition of $V_1=A\cup B$ with $|A|\le |B|\le |A|+1$ such that $\h[V_1]$ is $(\e/2)$-close to $\B_{|V_1|,k}(A, B)$ or $\overline \B_{|V_1|,k}(A, B)$. We add the vertices of $V\setminus V_1$ arbitrarily and evenly to $A$ and $B$, and get a partition of $V(\h)=A'\cup B'$ with $|A'|=\lfloor n/2\rfloor$, $A\subseteq A'$, and  $B\subseteq B'$.
Since $|V\setminus V_1|\le 4k\r^{14}n$, we conclude that $\h$ becomes a copy of $\B_{n,k}(A', B')$ or $\overline \B_{n,k}(A', B')$ after adding or deleting at most $\frac{\e}2 |V_1|^k+4k\r^{14}n\binom n{k-1}<\e n^k$ edges because $\r\le \e$. This means that $\h$ is $\e$-close to $\B_{n,k}$ or $\overline \B_{n,k}$, a contradiction. 

Furthermore, as $|V\setminus V_1|\le 4k\r^{14}n$, we have $\delta_d(\h[V_1])\ge \left(1/2-2\r\right) \binom{|V_1|-d}{k-d}$. We now apply Lemma \ref{lemR} on $\h[V_1]$ and get a  family $\R$ of order $(2\r)^{28}n$. 
Let $V_2 := V\setminus (V(\mathcal P_0)\cup V(\R))$, $n_2:=|V_2|$, and $\h_2: = \h[V_2]$. Note that $|V(\mathcal P_0)\cup V(\R)|\le 4k\r^{14}n +(3k/2)(2\r)^{28}n <\r^{13} n$ and thus $\delta_{d}(\h_2) \ge \left(1/2- 2\r\right) \binom{n_2-d}{k-d}$.
We now apply Lemma~\ref{lemP} to find a family of at most $p$ paths $\mathcal P_1, \mathcal P_2, \dots,\mathcal P_p$ covering all but at most $\a n_2$ vertices in $V_2$. For every $i\in [p]$, let $S_i$ and $T_i$ be two ends of $\mathcal P_i$. 
Due to Lemma \ref{lemR}, we can connect $S_i$ and $T_{i+1}$, $0\le i\le p$ (with $T_{p+1}:=T_0$), by disjoint $(3k/2)$-sets from $\R$ and get a $(k/2)$-cycle. This is possible because $p+1\le \r^{33}n\le \r^{32}|V_1|/3$.
At last, we use $\mathcal P_0$ to absorb all uncovered vertices in $V_2$ and unused vertices in $\R$.  
This is possible because the number of absorbed vertices is at most $\a n+(3k/2)|\R|\le 2k(2\r)^{28}n<k\r^{27}n/6$, and by our construction, this number is divisible by $k/2$.
\end{proof}

\medskip
It remains to prove the lemmas. 
We prove Lemmas~\ref{lemR} and~\ref{lemA} in Section~\ref{sec:4} via a Connecting Lemma, Lemma~\ref{lemC}, which itself is proved in Section~\ref{sec:C}.
In Section~\ref{sec:5} we introduce the weak regularity lemma and apply it to prove Lemma~\ref{lemP}.

\subsection{Proofs of Lemmas \ref{lemR} and \ref{lemA}}\label{sec:4}
Let us first state our connecting lemma and postpone its proof to Section~\ref{sec:C}.
\begin{lemma}[Connecting Lemma]\label{lemC}
For an even integer $k\ge 4$ and an integer $d$ with  $k/2\le d\le k-1$,
suppose $1/n\ll \r\ll\e\ll 1/k$.
Let $\h=(V,E)$ be a $k$-graph on $n$ vertices with $\delta_{d}(\h) \ge (1/2-\r) \binom{n-d}{k-d}$. 
If $\h$ is not $\e$-close to $\B_{n,k}$ or $\overline \B_{n,k}$, then there are at least $\r^4 n^{3k/2}$ connecting $(3k/2)$-sets for any two $(k/2)$-sets.
\end{lemma}

Now we can derive Lemmas~\ref{lemR} and~\ref{lemA} from Lemma~\ref{lemC}.
We use a concentration result from ~\cite{ABHKP16} for selecting connecting sets and absorbing paths.
Alternatively, we may use other well-known approaches, e.g., selecting sets uniformly at random and then removing the overlapping ones.

\begin{lemma}[Lemma~2.2 in~\cite{ABHKP16}]
  \label{lem:coupling} Let $\Omega$ be a finite probability space and
  let $\F_0\subseteq \dots\subseteq \F_n$ be partitions of~$\Omega$.  For each $i\in[n]$ let $Y_i$ be a Bernoulli random
  variable on $\Omega$ that is constant on each part of~$\F_i$, that
  is, let~$Y_i$ be $\F_i$-measurable.  Furthermore, let~$p_i$ be a
  real-valued random variable on~$\Omega$ 
  which is constant on each part of $\F_{i-1}$.
  Let~$x$ and~$\delta$ be  real numbers with $\delta\in(0,3/2)$, and
  let $X=Y_1+\dots+Y_n$.  If $\sum_{i=1}^n p_i\ge x$ holds almost
  surely and $\mathbb E [Y_i\mid\F_{i-1}]\ge p_i$ holds almost surely for
  all $i\in[n]$, then
  $\mathbb P\big(X<(1-\delta)x \big)<e^{-\delta^2 x/3}\,.$\qed
\end{lemma}

We would like to use Lemma~\ref{lem:coupling} to construct vertex-disjoint structures, that is, each time we select a vertex set (a connecting set or an absorbing path) uniformly at random from the ones 
disjoint from the previously chosen. 
For example, to construct the connecting $(3k/2)$-sets in Lemma~\ref{lemR}, let $\Omega$ be the collection of sequences $(S_1,\dots, S_n)$ of disjoint $(3k/2)$-sets, which are all possible outcomes of the sequential selection process.
Then the partitions $\F_0,\F_1,\dots, \F_n$ are defined as the `history' of the processes, namely, $\F_0=\Omega$ and for $i\in [n]$ $\F_i$ consists of collections of all sequences that share the same first $i$ terms of $(3k/2)$-sets.
Fix two $(k/2)$-sets $X_1$ and $X_2$, let $Y_i$ be the Bernoulli random variable on $\Omega$ such that it equals $1$ if and only if $(S_1,\dots, S_n)\in \Omega$ is such that $S_i$ is a connecting $(3k/2)$-set for $X_1$ and $X_2$ (clearly $Y_i$ is constant on each part of $\F_i$).
Similar setup can be used in our other applications, namely, in the proof of Lemma~\ref{lemA} and Claims~\ref{claim:absorbpath2} and~\ref{claim:absorbpath3}.
Finally, in all our applications, we will see that $p_i$ can be taken as a constant $p$ on $\Omega$, that is, we have $\mathbb E [Y_i\mid\F_{i-1}]\ge p$.

\begin{proof}[Proof of Lemma \ref{lemR}]
Suppose $1/n\ll \r\ll\e\ll 1/k$.
Let $\h$ be a $k$-graph on $n$ vertices such that $\delta_{d}(\h)\ge (\frac 12-\r) \binom{n-d}{k-d}$ and $\h$ is not $\e$-close to $\B_{n,k}$ or $\overline \B_{n,k}$. 
To find the family $\mathcal R$ we choose  $t= \r^{28} n$ disjoint connecting $(3k/2)$-sets of vertices $(S_1,\dots, S_t)$ and 
do so by sequentially selecting a uniformly random $(3k/2)$-set, which is connecting for some pair of $(k/2)$-sets and which is  disjoint from the previously chosen sets. 
For every two $(k/2)$-sets $X_1$ and $X_2$ and  every  $i\leq t$, 
let $\mathcal C_i$ 
be the collection of connecting sets 
for $X_1$ and $X_2$, which are disjoint from $S_1\cup\dots\cup S_{i-1}$.
By Lemma~\ref{lemC} the probability that $S_i$ is in $\mathcal C_i$ is at least $\r^4 - (3k/2)\r^{28} \ge \r^{4}/2$.
Thus, by Lemma~\ref{lem:coupling} 
with $\delta=1/3$ and $x=(\r^{4}/2)\r^{28}n$, with probability at least $1-e^{-\delta^2 x/3}$ the chosen family contains $(1-\delta)(\r^{4}/2)\r^{28}n = \r^{32}n/3$ connecting sets for each pair of $(k/2)$-sets.
Since $n^k e^{-\delta^2 x/3}<1$ the union bound implies that there exists a family $\R$ satisfying the property above for all pairs of $(k/2)$-sets simultaneously.
\end{proof}

Next we prove the Absorbing Lemma.
\begin{proof}[Proof of Lemma \ref{lemA}]
Suppose $1/n\ll \r\ll\e\ll 1/k$
and suppose $\h$ is a $k$-graph on $n$ vertices with $\delta_{d}(\h)\ge (\frac12 - \r)\binom {n-d}{k-d}$. 
By Proposition \ref{prop:seq}, we have that $\delta_{k/2}(\h)\ge (\frac12 - \r)\binom {n-k/2}{k/2}$.

 Given a set $X$ of $k/2$ vertices, an \emph{$X$-absorbing path~$\mathcal P$} is a $(k/2)$-path on $5k/2$ vertices such that there is a $(k/2)$-path on $V(\mathcal P)\cup X$, which has the same ends
  as~$\mathcal P$. 
 The core of the proof is the following claim showing that for any $(k/2)$-set $X$, there are many $X$-absorbing paths.

\begin{claim}\label{clm:9}
For any $(k/2)$-set $X$, there are at least $\r^{13}n^{5k/2}$ $X$-absorbing paths $\mathcal P$.
\end{claim}

\begin{proof}
For any $(k/2)$-set $X$ we construct the $X$-absorbing paths as follows. First choose a~$(k/2)$-set~$A$ such that $A\cup X\in E(\h)$ 
and note that there there are at least $\delta_{k/2}(\h)\ge (\frac12 - \r)\binom {n-k/2}{k/2}$ choices for~$A$. From $V(\h)\setminus (X\cup A)$
we choose two disjoint $(k/2)$-sets $B$ and $D$ such that $AB, DX\in E(\h)$ and such that $B$ and~$D$ have at least ${\r^4} n^{k/2}/2$ common neighbors $C$.
Note that  each of the choices $B,C, D$ yields a connecting $(3k/2)$-set for $A$ and $X$. Thus
there are at least ${\r^4} n^k/2$ choices for the pair $B, D$, as otherwise there are fewer than 
\[
\frac{\r^4}2 n^k\cdot n^{k/2} + n^k \cdot \frac{\r^4}2 n^{k/2} = {\r^4} n^{3k/2}
\]
connecting $(3k/2)$-sets for $A$ and $X$, which contradicts  Lemma \ref{lemC}. 

We pick $B, D$ as above and  pick  two disjoint common neighbors $C, E$ of $B, D$. 
Let $\mathcal P=ABCDE$ and note that $AXDCBE$ is also a $(k/2)$-path on $V(\mathcal P)\cup X$ with the same ends $A, E$ as $\mathcal P$. 
Moreover, the number of such $(5k/2)$-sets is at least
\[
\left(\frac12 - \r\right)\binom {n-k/2}{k/2}\cdot \frac{\r^4}2 n^k \cdot \binom{{\r^4} n^{k/2}/2}2 \ge \r^{13}n^{5k/2},
\]
as $\r$ is small enough.
\end{proof}

We choose a family of $\r^{14} n$ disjoint  $(5k/2)$-sets of vertices, 
doing so by sequentially selecting uniformly random $(5k/2)$-sets, which are absorbing for some $(k/2)$-set and which is  disjoint from the previously chosen ones. 
Note that for every $(k/2)$-set $X$ and in each step, by Claim~\ref{clm:9} the probability that the chosen $(5k/2)$-set is an $X$-absorbing path is at least $\r^{13} - (5k/2)\r^{14} \ge \r^{13}/2$.
Thus, by Lemma~\ref{lem:coupling} with $\delta=1/3$ and $x=(\r^{13}/2)\r^{14}n$ and the union bound, with probability at least $1-n^{k/2} e^{-\delta^2 x/3}>0$ the family $\mathcal F$ contains $(1-\delta)x = \r^{27}n/3$ $X$-absorbing paths for all $(k/2)$-sets $X$ simultaneously.
We take such a family and delete the $(5k/2)$-sets that are not absorbing paths for any $(k/2)$-set and connect the remaining $(5k/2)$-sets by Lemma \ref{lemC}. Since $3k/2$ vertices are used to connect each pair of $(5k/2)$-sets, we obtain the desired absorbing path which contains at most $\r^{14}n\cdot (5k/2+3k/2) =4k\r^{14}n$ vertices which can absorb at least $\r^{27}n/3$ $(k/2)$-sets, thus at least $k\r^{27}n/6$ vertices, proving the lemma. 
\end{proof}

\subsection{Proof of Lemma \ref{lemC}}\label{sec:C}
In this section we prove Lemma \ref{lemC}. Throughout this section we will use the following notation. Let $k\ge 4$ be an even integer. Given a $k$-graph $\h$, let $X=Y=\binom{V(\h)}{k/2}$. Set $N:=|X|=\binom n{k/2}$.

Given a $k$-graph $\h$, we define the bipartite graph $G(\h)$ as follows: $G(\h)$ has vertex classes $X$ and $Y$. Two vertices $x\in X$ and $y\in Y$ are adjacent in $G(\h)$ if and only if $x\cup y\in E(\h)$. When it is clear from the context, we will refer to $G(\h)$ as $G$.

Let $n, k\ge 4$ be positive integers with $k$ even. Denote by $B_{n,k}$ the bipartite graph with vertex classes $X$ and $Y$ both of sizes $N$ that satisfies the following properties: 
\begin{itemize}
\item
$X_1, X_2$ is a partition of $X$ such that $|X_1|=\lfloor N/2\rfloor $ and $|X_2|= \lceil N/2 \rceil$.
\item
$Y_1, Y_2$ is a partition of $Y$ such that $|Y_1|=\lfloor N/2\rfloor $ and $|Y_2|= \lceil N/2 \rceil$.
\item
$B_{n,k}[X_1, Y_1]$ and $B_{n,k}[X_2, Y_2]$ are complete bipartite graphs. Furthermore, there are no other edges in $B_{n,k}$.
\end{itemize}

We will use the following lemma from \cite{TrZh13}.

\begin{lemma}[Lemma 5.4, \cite{TrZh13}]\label{lem:TZ}
Given any $\e>0$ and even integer $k\ge 4$, there exist $\beta>0$ and $n_0\in \mathbb N$ such that the following holds. Suppose that $\h$ is a $k$-uniform hypergraph on $n\ge n_0$ vertices. Suppose further that $G:=G(\h)$ satisfies $G=B_{n,k}\pm \beta N^2$. Then $\h$ is $\e$-close to $\B_{n,k}$ or $\overline \B_{n,k}$.
\end{lemma}

The next claim shows that under our degree condition, if two $(k/2)$-sets have many connecting $(k/2)$-sets, then they have many connecting $(3k/2)$-sets.

\begin{claim}\label{clm:1-3}
Suppose $1/n\ll \r\ll 1/k$.
Let $\h$ be a $k$-graph with $\delta_{k/2}(\h)\ge (\frac12 -\r)\binom{n-k/2}{k/2}$.
If two $(k/2)$-sets $x, y\in V(\h)$ have at least $\r \binom{n}{k/2}$ connecting $(k/2)$-sets, then they have at least $\r^4 n^{3k/2}$ connecting $(3k/2)$-sets. 
\end{claim}

\begin{proof}
First note that by $\delta_{k/2}(\h)\ge (\frac12 -\r)\binom{n-k/2}{k/2}$, for any three $(k/2)$-sets $a, b, c$, at least two of them have at least $(\frac16-2\r)\binom{n-k/2}{k/2}$ common neighbors in $\h$. 
Indeed, for any set $S$ in $V(\h)$, let $N_\h(S)$ be the collection of $(k-|S|)$-sets $T$ in $V(\h)\setminus S$ such that $T\cup S\in E(\h)$ and
assume that $|N_\h(a)\cap N_\h(b)|<(\frac16 - 2\r)\binom{n-k/2}{k/2}$, thus $|N_\h(a)\cup N_\h(b)|>\frac56\binom{n-k/2}{k/2}$. 
Together with the minimum degree condition, this implies that $|N_\h(c)\cap N_\h(a)|>(\frac16 - \r)\binom{n-k/2}{k/2}$ or $|N_\h(c)\cap N_\h(b)|>(\frac16 - \r)\binom{n-k/2}{k/2}$. 

To prove the claim, let $Z$ be the set of connecting $(k/2)$-sets for $x, y$, then $x\cup z, y\cup z\in E(\h)$ for any $z\in Z$. 
By the discussion in the previous paragraph, fix any three sets $z_1, z_2, z_3\in Z$, there exists one pair of them such that they have at least $(\frac16-2\r)\binom{n-k/2}{k/2}$ common neighbors in $\h$, so at least $\frac17\binom{n-k/2}{k/2}$ common neighbors in $V(\h)\setminus\{x\cup y\cup z_1\cup z_2\cup z_3\}$. 
We count the number of connecting $(3k/2)$-sets by taking the sum over all triples of sets in $\binom Z3$. 
Since every pair in $Z$ can be counted at most $|Z|-2$ times, we obtain
\[
\frac1{|Z|-2}\binom{|Z|}3\cdot \frac17\binom{n-k/2}{k/2} >\frac{|Z|^2}{49}\binom{n-k/2}{k/2} > \r^3 n^{3k/2}
\]
ordered multisets of $3k/2$ elements, where the last inequality is because $\r$ is small.
Since each $(3k/2)$-set can be counted at most $(3k/2)!$ times and a multiset with repeated elements contributes $O(n^{3k/2-1})$ to the quantity above, we obtain at least
\[
({(3k/2)!})^{-1}\r^3 n^{3k/2} -O(n^{3k/2-1}) \ge  \r^4 n^{3k/2}
\]
connecting $(3k/2)$-sets for $x$ and $y$.
\end{proof}

Now we are ready to prove Lemma \ref{lemC}.
For a given graph $G$ and two disjoint vertex subsets $A, B$, let $e_G(A, B)$ denote the number of edges of $G$ with one end in $A$ and one end in $B$ and let $d_G(A, B)=\frac{e_G(A, B)}{|A||B|}$.
The subscript will be omitted when the graph $G$ is clear from context.

\begin{proof}[Proof of Lemma \ref{lemC}]

Suppose $1/n\ll \r\ll\e\ll 1/k$.
By Proposition \ref{prop:seq}, it suffices to prove the lemma for $d=k/2$. Let $\h$ be a $k$-graph on $n$ vertices such that $\delta_{k/2}(\h)\ge (\frac12 -\r)\binom{n-k/2}{k/2}\ge \left(\frac12 - 2\r\right)\binom{n}{k/2}$ and $\h$ is not $\e$-close to $\B_{n,k}$ or $\overline \B_{n,k}$.

Now assume to the contrary, that there are two $(k/2)$-sets $x$ and $y$ with fewer than $\r^4 n^{3k/2}$ connecting $(3k/2)$-sets. 
Thus, by Claim \ref{clm:1-3}, there are fewer than $\r \binom{n}{k/2}$ connecting $(k/2)$-sets for $x$ and $y$. 

Consider $G:=G(\h)$ and we have $\delta(G)\ge (\frac12-2\r)N$. 
Consider $x, y\in X$ (note that $x,y$ also exist in $Y$). By our assumption, $|N_G(x)\cap N_G(y)|< \r N$. 
Let $Y_1\subset Y$ of size exactly $\lfloor N/2\rfloor$ such that $Y_1$ maximizes $|Y_1\cap N_G(x)|$ and thus minimizes $|Y_1\cap (N_G(y)\setminus N_G(x))|$. 
Let $Y_2:=Y\setminus Y_1$. By definition and $|N_G(x)\cap N_G(y)|< \r N$, it is easy to see that
\begin{equation}\label{eq:Y}
|Y_1\setminus N_G(x)|\le 2\r N, \quad |Y_2\setminus N_G(y)|\le 3\r N.
\end{equation}
Let $X'=\{z\in X\setminus\{x,y\}\colon |N_G(x)\cap N_G(z)|\ge \r N \text{ and } |N_G(y)\cap N_G(z)|\ge \r N\}$. We claim that $|X'|< \r N$. Indeed, otherwise, by greedily picking $z\in X'$, $a\in (N_G(x)\cap N_G(z))\setminus\{y\}$ and $b\in (N_G(y)\cap N_G(z))\setminus \{x,a\}$, we get at least $\r N(\r N-1) (\r N-2)> \r^4 n^{3k/2}$ connecting $(3k/2)$-sets $\{a, z, b\}$ for $x$ and $y$ in $\h$ (because $\r$ is small), a contradiction. Thus, since $\delta(G)\ge (\frac12 - 2\r)N$, for all but at most $\r N$ vertices $z\in X$, either 
\begin{enumerate}
\item[i)]
$|N_G(x)\cap N_G(z)|\ge \r N$ and $|N_G(y)\cap N_G(z)|< \r N$, or
\item[ii)]
$|N_G(x)\cap N_G(z)|< \r N$ and $|N_G(y)\cap N_G(z)|\ge \r N$.
\end{enumerate}
Indeed, since $| N_G(x) \cup N_G(y) |\ge (1- 5\r)N$, any vertex $z$ not in $X'$ and not satisfying i) or ii) satisfies that $|N_G(x)\cap N_G(z)|< \r N$ and $|N_G(y)\cap N_G(z)|< \r N$ and thus $\deg_G(z)\le \r N + \r N + 5\r N < (\frac12-2\r)N$, a contradiction.
Let $X_1$ be the set of vertices in $X$ satisfying property i) and $X_2$ be the set of vertices in $X$ satisfying property ii). Clearly, $X_1\cap X_2=\emptyset$. By definition and~\eqref{eq:Y}, for any $z\in X_1$,
\[
|N_G(z)\cap Y_2|\le |N_G(y)\cap N_G(z)| + |Y_2\setminus N_G(y)|< \r N + 3\r N = 4\r N. 
\]
Then by $\delta(G)\ge (\frac12-2\r)N$, we get $|N_G(z)\cap Y_1|> \left(\frac12 - 6\r\right)N$. Similarly, for any $z\in X_2$, we have $|N_G(z)\cap Y_1|< 3\r N$ and $|N_G(z)\cap Y_2|> \left(\frac12 - 5\r\right)N$. 
Together with $|Y_1|=\lfloor N/2\rfloor $ and $|Y_2|= \lceil N/2 \rceil$, we get
\begin{align}
&d(X_1, Y_2)< 8\r, \quad d(X_1,Y_1)\ge 1 - 12\r, \nonumber \\
&d(X_2, Y_1)< 6\r, \quad d(X_2, Y_2)\ge 1 - 10\r.   \label{eq:den}
\end{align}

We also claim that $|X_1|\ge (\frac12 - 9\r)N$ and $|X_2|\ge (\frac12 - 11\r)N$. 
Indeed, if $|X_1|< (\frac12 - 9\r)N$, by summing up the degrees of vertices in $Y_1$ and~\eqref{eq:den}, we get $\sum_{v\in Y_1}\deg_G (v)< |X_1|  |Y_1| + |X'| |Y_1| + 6\r |X_2| |Y_1| $.
By averaging, there is a vertex $v\in Y_1$ such that
\[
\deg_G (v)< \left(\frac12 - 9\r \right)N + \r N + 6\r N= \left( \frac12-2\r \right)N,
\]
a contradiction.
Similar calculations show that $|X_2|\ge (\frac12 - 11\r)N$. 
In summary, we get
\begin{equation}\label{eq:X12}
\left(\frac12 - 9\r \right)N \le |X_1|\le \left(\frac12 + 11\r \right)N \quad \text{and} \quad \left(\frac12 - 11\r \right)N \le |X_2| \le \left(\frac12 + 9\r \right)N.
\end{equation}
Finally, let $X_1'\subset X$ of size exactly $\lfloor N/2 \rfloor$ such that $X_1'$ maximizes $|X_1'\cap X_1|$ and minimizes $|X_1'\cap X_2|$.
Let $X_2':=X\setminus X_1'$. 
Thus we get a partiton of $X=X_1'\cup X_2'$, $Y=Y_1\cup Y_2$, where for $i=1,2$, $X_i'$ plays the role of $X_i$ as in the definition of $B_{n,k}$. 
We claim that $G=B_{n,k}\pm 30\r N^2$. 
Note that if this is true, by Lemma \ref{lem:TZ} with $\beta=30\r$, we get that $\h$ is $\e$-close to $\B_{n,k}$ or $\bar \B_{n,k}$, a contradiction.
This contradiction will complete the proof.

Indeed, if $|X_1|\ge N/2$, we have $X_1'\subseteq X_1$, so by the property of the vertices in $X_1$, $e(X_1', Y_2)< 4\r N\cdot |X_1'|=8\r (N/2)^2$, $e(X_1', Y_1)>(\frac12 - 6\r) N\cdot |X_1'|=(1 - 12\r)(N/2)^2$. 
Otherwise, $|X_1|< N/2$, implying that $X_1\subseteq X_1'$.
By \eqref{eq:den} and \eqref{eq:X12}, we infer that
\begin{align*}
e(X_1', Y_2)&< 8\r |X_1||Y_2|+|X_1'\setminus X_1||Y_2|< \left(8\r N/2 + 9\r N \right) |Y_2|=26\r (N/2)^2, \text{ and} \\
e(X_1', Y_1)&\ge e(X_1, Y_1) \ge  (1 - 12\r) \left(\frac12 - 9\r \right)N \left\lfloor \frac N2 \right\rfloor> (1-30\r) (N/2)^2.
\end{align*}
In both cases, we have $e(X_1', Y_2)< 26\r (N/2)^2$ and $e(X_1', Y_1) >(1-30\r) (N/2)^2$. 
Similarly, for $X_2'$, if $|X_2|\ge N/2$, then we have $X_2'\subseteq X_2$, so by the property of the vertices in $X_2$,  $e(X_2', Y_1)<3\r N\cdot |X_2'|=6\r (N/2)^2$ and $e(X_2', Y_2)>(\frac12 - 5\r)N |X_2'|=(1 - 10\r)(N/2)^2$. 
Otherwise, $|X_2|< N/2$, implying that $X_2\subseteq X_2'$.
By \eqref{eq:den} and \eqref{eq:X12}, we infer that
\begin{align*}
e(X_2', Y_1)&< 6\r |X_2||Y_1|+|X_2'\setminus X_2||Y_1|< \left(6\r N/2 + 11\r N \right) |Y_1|=28\r (N/2)^2, \text{ and} \\
e(X_2', Y_2)&\ge e(X_2, Y_2) \ge  (1 - 10\r) \left(\frac12 - 11\r \right)N \frac N2> (1-32\r) (N/2)^2.
\end{align*}
In both cases, we have $e(X_2', Y_1)< 28\r (N/2)^2$ and $e(X_2', Y_2)>(1-32\r) (N/2)^2$. 

In summary, for the partition $[X_1'\cup X_2', Y_1\cup Y_2]$ of $G$, we have, 
\[
e(X_1', Y_2), e(X_2', Y_1)< 28\r (N/2)^2 \text{ and } e(X_1', Y_1), e(X_2', Y_2)>(1-32\r) (N/2)^2.
\]
Thus, we conclude that $G=B_{n,k}\pm 30\r N^2$ because $(32\r +32\r +28\r +28\r) (N/2)^2=30\r N^2$. 
\end{proof}

\subsection{Proof of Lemma \ref{lemP}}\label{sec:5}

We 
follow the approach from \cite{HS}, which uses the weak regularity lemma for hypergraphs, a straightforward extension of Szemer\'edi's regularity lemma for graphs~\cite{Sze}. 

Let $\h = (V, E)$ be a $k$-graph and let $A_1, \dots, A_k$ be mutually disjoint non-empty subsets of $V$. We define $e(A_1, \dots, A_k)$ to be the number of edges with one vertex in each $A_i$, $i\in [k]$, and the density of $\h$ with respect to ($A_1, \dots, A_k$) as
\[
d(A_1,\dots, A_k) = \frac{e(A_1, \dots, A_k)}{|A_1| \cdots|A_k|}.
\]
Given $\e, d\ge 0$, a $k$-tuple ($V_1, \dots, V_k$) of mutually disjoint subsets $V_1, \dots, V_k\subset V$ is \emph{$(\e,d)$-regular} if
\[
|d(A_1, \dots, A_k) - d|\le \e
\]
for all $k$-tuples of subsets of $A_i\subset V_i$, $i\in [k]$, satisfying $|A_i|\ge \e |V_i|$. We say ($V_1, \dots, V_k$) is $\e$-regular if it is $(\e, d)$-regular for some $d\ge 0$. It is immediate from the definition that in an $(\e, d)$-regular $k$-tuple ($V_1, \dots, V_k$), if $V_i'\subset V_i$ has size $|V_i'| \ge c|V_i|$ for some $c\ge \e$, then ($V_1', \dots, V_k'$) is $(\e/c, d)$-regular.

\begin{theorem}\label{thmReg}
For all $t_0\ge 0$ and $\e>0$, there exist $T_0 = T_0(t_0, \e)$ and $n_0 = n_0(t_0,\e)$ so that for every $k$-graph $\h = (V, E)$ on $n>n_0$ vertices, there exists a partition $V = V_0 \dot \cup V_1 \dot \cup \cdots \dot\cup V_t$ such that
\begin{enumerate}
\item[(i)] $t_0\le t\le T_0$,
\item[(ii)] $|V_1| = |V_2| = \cdots = |V_t|$ and $|V_0|\le \e n$,
\item[(iii)] for all but at most $\e \binom tk$ sets $\{i_1,\dots, i_k\}\in \binom{[t]}{k}$, the $k$-tuple $(V_{i_1}, \dots, V_{i_k})$ is $\e$-regular.
\end{enumerate}
\end{theorem}

A partition as given in Theorem \ref{thmReg} is called an $(\e,t)$-regular partition of $\h$. For an $(\e, t)$-regular partition of $\h$ and $d\ge 0$ we refer to $\Q = (V_i)_{i\in [t]}$ as the family of \emph{clusters} and define the \emph{cluster hypergraph} $\K = \K(\e,d,\Q)$ with vertex set $[t]$ and $\{i_1,\dots,i_k\}\in \binom{[t]}{k}$ is an edge if and only if $(V_{i_1}, \dots, V_{i_k})$ is $\e$-regular and $d(V_{i_1}, \dots, V_{i_k})\ge d$.

The following corollary shows that the cluster hypergraph inherits the minimum codegree of the original hypergraph.
The proof is standard and very similar to that of \cite[Proposition 16]{HS} so we omit the proof.

\begin{corollary}\label{prop16}
For $c, \e, d>0$, an even integer $k\ge 3$, and an integer $t_0\ge 2k^2/d$, there exist $T_0$ and $n_0$ such that the following holds. Given a $k$-graph $\h$ on $n>n_0$ vertices with $\delta_{k/2}(\h)\ge c \binom{n-k/2}{k/2}$, there exists an $\e$-regular partition $\Q = (V_i)_{i\in [t]}$,  with $t_0\le t\le T_0$. Furthermore, let $\K = \K(\e,d/2,\Q)$ be the cluster hypergraph of $\h$. Then the number of $(k/2)$-sets $S\in \binom{[t]}{k/2}$ violating $\deg_{\K}(S)\ge (c - \sqrt{\e} - d)\binom{t-k/2}{k/2}$ is at most $\sqrt{\e} \binom{t}{k/2}$.
\end{corollary}

We use the following proposition from \cite[Claim 4.1]{RRS08}.

\begin{proposition}\label{prop12}
Given $d>0$ and $k\ge 2$, every $k$-partite $k$-graph $\h$ with at most $m$ vertices in each part and with at least $dm^k$ edges contains a $(k/2)$-path on at least $dm$ vertices.
\end{proposition}

We want to use Proposition \ref{prop12} to cover an $(\e,d)$-regular tuple $(V_1,\dots, V_k)$ by $(k/2)$-paths. 
Note that a $k$-partite $(k/2)$-path of odd length $t$ has $(t+1)/2$ vertices in each cluster.

\begin{lemma}\label{lem:path}
Fix an even integer $k\ge 4$ and $\e, d>0$ such that $d>2\e$. Let $m>\frac{k}{\e(d-\e)}$. Suppose $\mathcal V = (V_1, V_2,\dots, V_k)$ is an ($\e,d$)-regular $k$-tuple with $|V_i| = m$ for $i \in [k]$. Then there is a family consisting of $\frac{k}{(d-2\e)\e}$ pairwise vertex-disjoint $(k/2)$-paths which cover all but at most $k\e m$ vertices of $\mathcal V$.
\end{lemma}

\begin{proof}
We greedily find $(k/2)$-paths of odd length by Proposition \ref{prop12} in $\mathcal V$ until every cluster has less than $\e m$ vertices uncovered. Assume that every cluster has $m'\ge \e m$ vertices uncovered. 
By regularity, the remaining hypergraph has at least $(d-\e)(m')^k$ edges. 
We apply Proposition~\ref{prop12} and get a $(k/2)$-path of odd length covering at least $(d-\e)m'-k/2\ge (d-2\e)\e m$ vertices (we discard one $(k/2)$-set if needed). 
Thus, the number of paths is at most $km/((d-2\e)\e m)= \frac{k}{(d-2\e)\e}$.
\end{proof}

We will find an almost perfect matching in the cluster hypergraph. 
In \cite{MaRu}, Markstr\"{o}m and Ruci\'nski stated the following theorem for $1\le d<k/2$ and assumed a minimum $d$-degree condition for all $d$-sets. In fact, their proof works for all $d$ with $1\le d\le k-2$ and can be easily adapted to prove the following theorem, in which a small collection of $d$-sets are allowed to have degree zero.

\begin{theorem}\label{thm:MaRu}
For each integer $k\ge 3$, $1\le d\le k-2$ and every $0<\r<1/4$, $\e>0$ the following holds for sufficiently large $n$. Suppose that $\h$ is a $k$-graph on $n$ vertices such that for all but at most $\e \binom{n}{d}$ $d$-sets $S$,
\[
\deg(S)\ge \left( \frac{k-d}{k} - \frac{1}{k^{k-d}}+\r \right) \binom{n-d}{k-d}.
\]
Then $\h$ contains a matching that covers all but at most $2 \e^{1/k}n$ vertices.
\end{theorem}

\begin{proof}
Let $\M$ be a largest matching in $\h$. 
Assume to the contrary that $n-|V(\M)|\ge 2 \e^{1/k}n$. Let $X:=V(\h)\setminus V(\M)$ and $m:=|\M|$. 
We call a $d$-set $S\in \binom Xd$ bad if $\deg(S)< \left( \frac{k-d}{k} - \frac{1}{k^{k-d}}+\r \right) \binom{n-d}{k-d}$. So the number of bad $d$-sets in $\h$ is at most $\e \binom{n}{d}$.

For every $S\in \binom Xd$ and any submatching $\M'$ of $\M$, denote by $L_S(\M')$ the $(k-d)$-graph consisting of all $(k-d)$-sets $T\subseteq V(\M')$ such that $S\cup T\in E(\h)$ and $|T\cap e|\le 1$ for every edge $e\in \M'$. 
Note that for a fixed set $S\subseteq X$ we have $L_S(\M)=\bigcup_{E\in \binom{\M}{k-d}}L_S(E)$, where the hypergraphs $L_S(E)$ are pairwise edge-disjoint.

For every $S\in \binom Xd$, we break the family $\binom{\M}{k-d}$ consisting of the sets $E=\{e_1,\dots, e_{k-d}\}$, where $e_i\in \M$, into three parts, according to the properties of the link $L_S(E)$. Namely, we write $\binom{\M}{k-d}=P(S)\cup A(S)\cup B(S)$, where:
\begin{itemize}
\item $P(S):=\left\{E\in \binom{\M}{k-d}: L_S(E)\text{ has a matching of size }k-d+1 \right\}$.
\item $A(S):=\left\{E\in \binom{\M}{k-d}:|L_S(E)|\le (k-d)k^{k-d-1}-1 \right\}$.
\item $B(S):=\left\{E\in \binom{\M}{k-d}\setminus P(S): |L_S(E)|=(k-d)k^{k-d-1} \right\}$.
\end{itemize}

We omit the proofs of these two facts because the minimum degree condition is not involved in their proofs \cite[Facts 2 and 3]{MaRu}.

\begin{itemize}
\item[\textit{Fact 1.}] For at most $\r\binom{|X|}{d}$ sets $S\in\binom Xd$ we have $|P(S)|>\frac13 \r \binom{m}{k-d}$. 
\item[\textit{Fact 2.}] For at most $\r\binom{|X|}{d}$ sets $S\in\binom Xd$ we have $|B(S)|>\frac13 \r \binom{m}{k-d}$. 
\end{itemize}

By these two facts, for all but at most $2\r\binom{|X|}{d}$ $d$-sets $S\in \binom Xd$, we have $|P(S)|+ |B(S)|\le\frac23 \r \binom{m}{k-d}$, thus,
\begin{align*}
|L_S(\M)|=\sum_{E\in \binom{\M}{k-d}}L_S(E)&\le k^{k-d}(|P(S)|+|B(S)|)+((k-d)k^{k-d-1}-1)|A(S)|\\
&\le \left( \frac{2\r}{3}k^{k-d}+(k-d)k^{k-d-1}-1 \right) \binom{m}{k-d}\\
&\le \left( \frac{2\r}{3}+\frac{k-d}{k}-\frac{1}{k^{k-d}} \right) \binom{n}{k-d} \quad \text{ by }m\le \frac nk
\end{align*}
where we used the trivial bounds $|L_S(E)|\le k^{k-d}$ and $|A(S)|\le \binom m{k-d}$. 
Observe that given such an $S$, the number of $(k-d)$-sets $T$ such that $S\cup T\in E(\h)$ and $T\notin L_S(\M)$ is $o(n^{k-d})$. Hence, we have
\begin{align*}
\deg(S)&=|L_S(\M)|+ o(n^{k-d})\le \left( \frac{2\r}{3}+\frac{k-d}{k}-\frac{1}{k^{k-d}}+o(1) \right) \binom{n}{k-d}\\
&<\left( \frac{k-d}{k} - \frac{1}{k^{k-d}}+\r \right) \binom{n-d}{k-d},
\end{align*}
which means $S$ is bad. Since $|X|\ge 2 \e^{1/k}n$, this implies that the number of bad $d$-sets in $X$ is at least
\[
\binom{|X|}{d} - 2\r\binom{|X|}{d}>\frac12\binom{|X|}{d}\ge \frac12\binom{2 \e^{1/k}n}{d}>\e^{d/k}\binom{n}{d}>\e\binom nd,
\]
a contradiction.
\end{proof}

Now we are ready to prove Lemma \ref{lemP}.

\begin{proof}[Proof of Lemma \ref{lemP}]
Let $k, d$ be integers such that $k\in 2\mathbb N$ and $k/2\le d\le k-1$.
Suppose $1/n \ll 1/p\ll 1/T_0 \ll 1/t_0\ll \e \ll \a \ll \r \ll 1/k$. 

It suffices to prove the lemma for the case $d=k/2$. 
Suppose $\h$ is a $k$-graph on $n$ vertices and $\delta_{k/2}(\h)\ge (\frac12 - \r)\binom{n-k/2}{k/2}$. 
We apply Corollary~\ref{prop16} with parameters $\frac12 - \r$, $\e$, $2\r$ and $t_0$ obtaining an $(\e, t)$-regular partition $\Q=(V_i)_{i\in [t]}$ with $t_0\le t\le T_0$ and the cluster hypergraph $\K=\K(\e, \r, \Q)$ with vertex set $[t]$. 
Let $m\ge \frac{(1-\e)n}{t}$ be the size of each cluster $V_i$, $i\in [t]$. 
By Corollary \ref{prop16}, for all but at most $\sqrt{\e}\binom{t}{k/2}$ $(k/2)$-sets $S$, 
\[
\deg_{\K}(S)\ge \left(\frac12 - \r - \sqrt{\e} - 2\r \right)\binom{t-k/2}{k/2}\ge\left(\frac12 - 4\r \right)\binom{t-k/2}{k/2}.
\]
Note that we have $\frac12 - 4\r>\frac{k/2}{k} - \frac{1}{k^{k/2}}+\r$ because $\r$ is small. 
Thus by Theorem \ref{thm:MaRu}, $\K$ contains a matching $\M$ covering all but at most $2 \e^{1/k}t$ vertices. 
For each edge $\{i_1,\dots,i_k\}\in \M$, the corresponding clusters $(V_{i_1}, \dots, V_{i_k})$ is $(\e, \r')$-regular for some $\r'\ge \r$. 
Thus we can apply Lemma \ref{lem:path} on $(V_{i_1}, \dots, V_{i_k})$ and get a family of at most $\frac{k}{(\r-2\e)\e}$ $(k/2)$-paths leaving at most $k\e m$ vertices uncovered. 
We do this for each edge in $\M$ and get at most $\frac{t}{k}\cdot \frac{k}{(\r-2\e)\e}\le p$ $(k/2)$-paths, which leaves at most
\[
|V_0|+k\e m\cdot \frac tk+2 \e^{1/k}t \cdot m\le \e n + \e n + 2 \e^{1/k}n<3\e^{1/k}n\le \a n
\]
vertices uncovered in $\h$.
\end{proof}

\section{The extremal case - proof of Theorem \ref{lemE} }
\label{sec:8}

This section is devoted to the proof of Theorem \ref{lemE}.
For two  $k$-graphs $\h, \h'$  on the same vertex set $V$, let $\h'\setminus \h:=(V, E(\h')\setminus E(\h))$. Suppose that $0\le \a\le 1$ and $|V|=n$.  We call a set $S\subset V$ in $\h$ 
\begin{itemize}
\item \emph{$\a$-good} with respect to $\h'$  if $\deg_{\h'\setminus \h}(S) \le \a n^{k-|S|}$;
\item \emph{$\a$-bad} with respect to $\h'$  if $\deg_{\h'\cap \h}(S) \le \a n^{k-|S|}$;
\item \emph{$\a$-medium}  with respect to $\h'$ otherwise.
\end{itemize}
Let $\h$  be given as in Theorem~\ref{lemE}. In particular, 
\begin{align} \label{eq:md}
\delta_{d}(\h) > \overline\delta(n,k,d) = \max \{\delta_d(\F):  \F\in \h_{\ext}(n, k) \}.
\end{align}
When $\h$ is close to $\overline\B_{n,k}(A, B)$ (respectively, $\B_{n,k}(A, B)$), we call the even (respectively, odd) edges of $\h$ \emph{majority edges} and odd (respectively, even) edges \emph{minority edges}. 
By \eqref{eq:md} and Proposition~\ref{prop:seq}, we have $\delta_1(\h)\ge (1/2 - \e) \binom{n-1}{k-1}$. 

\subsection{Proof overview}\label{sec:overview}
Before delving into the details, we give an overview of the proof. It consists  of  three steps and we will also introduce some auxiliary results for the last step, whose proofs we defer to Section~\ref{sec:deferred}.
\begin{itemize}
\item[\textit{Step 0.}] Move all  vertices that are not $\a$-good with respect to $\B_{n,k}(A, B)$ (respectively, $\overline\B_{n,k}(A, B)$) and that are contained in more minority edges than majority edges to the other part to obtain a partition 
$A_1\cup B_1$ of $V$ with $|A_1|, |B_1|\approx n/2$ such that almost all vertices are good and no vertex is bad with respect to $\B_{n,k}(A_1, B_1)$ (or $\overline\B_{n,k}(A_1, B_1)$).
Moreover, every vertex is in $\frac15 \binom{n-1}{k-1}$ majority edges and almost all $(k/2)$-sets are good.
\item[\textit{Step 1.}] Build a constant size path $\pp_{\mathrm b}$ which  breaks the parity barriers discussed in Section 1.3 and extend it to a short path $\pp$ which contains all the medium vertices.
As mentioned in Section 1.3 breaking the parity barrier is a crucial part of the proof and the construction of $\pp_{\mathrm b}$  is split into several cases depending on whether $n\in \frac k2\mathbb N\setminus k \mathbb N$ or not, and whether $\h$ is $\e$-close to $\overline{\B}_{n,k}$ or ${\B}_{n,k}$.
In each case, we apply Lemma~\ref{prop:e123} or Lemma~\ref{thm:deg_bridge} to break the parity barriers.

\item[\textit{Step 2.}] Let  $L$ and $S$ denote the ends of the path $\pp$ (which are sets of size $k/2$) from Step 1 and by possibly extending $\pp$, using more vertices from the larger one of $A_1$ and $B_1$, we make sure
that $|A_1\setminus V(\pp)|=|B_1\setminus V(\pp)|$. Our goal is to pick an edge $S'L'$ from $V\setminus V(\pp)$ 
and find suitable paths $\pp_1$ with ends $S, S'$ and $\pp_2$ with ends~$L,L'$ such that $L~\pp~S~\pp_1~S'~L'~\pp_2~L$ forms a Hamilton $(k/2)$-cycle of $\h$.
Below we give some details on how to obtain the paths $\pp_1$ and $\pp_2$.

We call a set  $T\subseteq V(\h)$ an \emph{$(i,j)$-set (wrt.~$A_1$ and $B_1$)}  if $|T\cap A_1|=i$ and $|T\cap B_1|=j$, and an edge an $(i,j)$-edge if it is an $(i,j)$-set.
Given integers $0\le r\le k$ and two  sets $A$ and~$B$, let $\K^k(A)$ be the complete $k$-graph on~$A$ and  let $\K_r^k(A, B)$ 
be the $k$-graph on $A\cup B$ whose edges are all $k$-sets intersecting~$A$ in precisely $r$ vertices. 
We pick an edge $S'L'$ from $V\setminus V(\pp)$ and for simplicity we assume in the following that $S,S'\subset A_1$ and $L,L'\subset B_1$.

First suppose that $\h$ is close to ${\B}_{n,k}$ and
note that in this case almost all $(k-1,1)$-sets and $(1,k-1)$-sets are edges of $\h$.
We use them to build two long paths, one with ends $L, L'$ and consisting of $(k/2, 0)$-sets and $(k/2-1,1)$-sets alternately and  
the other with ends $S, S'$ and consisting of $(0, k/2)$-sets and $(1, k/2-1)$-sets alternately.
To achieve this, we essentially split $A_1\setminus V(\pp)$ and $B_1\setminus V(\pp)$ each into two parts with ratio $1:(k-1)$ and apply the following lemma twice to obtain $\pp_1$ and $\pp_2$.

\begin{lemma}\label{lem:finish2}
Given an even integer $k\ge 4$, suppose $1/t\ll \a_0\ll 1/k$. 
Suppose that $\h$ is a $k$-graph on $V=X\dot\cup Y$ such that $|X|=t$, $|Y|=t(k-1)+k/2$, and every vertex of $\h$ is $\a_0$-good with respect to $\K_1^k(X, Y)$. 
Then, given any two disjoint $(k/2)$-sets $L_0, L_1\subset Y$, which are ${\a_0}$-good   with respect to $\K_1^k(X, Y)$, there is a Hamilton $(k/2)$-path in $\h$ with ends $L_0$ and $L_1$. 
\end{lemma}

Now suppose that $\h$ is close to $\overline{\B}_{n,k}$.
In this case the $(k/2)$-sets we use for the Hamilton cycle  must have the same parity as the ends of $\pp$ and
we will have to deal with four cases depending on the parity of $k/2$ and that of the ends of $\pp$.
In order to illustrate the main ideas, we  elaborate  on the case when both ends of $\pp$ are even.

If $k/2$ is even then
almost all $(k,0)$-sets and $(0,k)$-sets are  edges of $\h$ and we find the two paths~$\pp_1$ and $\pp_2$ by using the following lemma.
\begin{lemma}\label{lem:finish}
Given an even integer $k\ge 4$, suppose $1/t\ll \a_0\ll 1/k$, where $t$ is an integer.
Suppose that $\h$ is a $k$-graph on $Y$ of order $|Y| = kt/2$ such that every vertex is $\a_0$-good with respect to $\K^k(Y)$. 
Then, given any two disjoint $(k/2)$-sets  $L_0, L_1 \subset Y$, which are ${\a_0}$-good   with respect to $\K^k(Y)$, there is a Hamilton $(k/2)$-path in $\h$ with ends $L_0$ and $L_1$. 
\end{lemma}

Next, assume that $k/2$ is odd.
As a $(k,0)$-set can only be split into two odd $(k/2)$-sets yet $\h$ is close to $\overline{\B}_{n,k}$, these sets are not useful.
Instead we will use the $(k-2,2)$-sets and $(0,k)$-sets to construct one path with ends $L, L'$ and consisting of $(k/2-1, 1)$-sets, the other with ends $S, S'$ and consisting of $(0, k/2)$-sets.
Thus, we essentially split $A_1\setminus V(\pp)$ and $B_1\setminus V(\pp)$ each into two parts with ratio $2:(k-2)$ and apply the following lemma twice to obtain $\pp_1$ and $\pp_2$.
\begin{lemma}\label{lem:finish1}
Given an even integer $k\ge 4$, suppose $1/t\ll \a_0\ll 1/k$. 
Suppose that $\h$ is an $n$-vertex $k$-graph on $V=X\dot\cup Y$ such that $|X|=t$, $|Y|=(k/2-1)t$, and that every vertex of $\h$ is $\a_0$-good with respect to $\K_2^k(X, Y)$. 
Then, given any two disjoint $(k/2)$-sets $L_0, L_1$, $|L_i\cap X|=1$, $i\in\{0,1\}$,  which are ${\a_0}$-good  with respect to $\K_2^k(X, Y)$, there is a Hamilton $(k/2)$-path in $\h$ with ends $L_0$ and~$L_1$.
\end{lemma}
\end{itemize}
As mentioned above, we postpone the proofs of  Lemmas~\ref{lem:finish2}--\ref{lem:finish1} to Section~\ref{sec:deferred} and first continue with the details of Step 0 and Step 1.

Throughout the rest of the paper let $\e_0:=k^{\frac12}\e^{1-\frac1{2k}}$, $\e_1:=k^{\frac12}\e^{\frac1{2k}}$, $\e_0':=k \e^{1-\frac1k}$ and $\e':=k^{\frac12}\e^{\frac14}$. 

\subsection{Step 0 - Finding a suitable partition.} \label{sec:step0}
Let $A\cup B$ be a partition of $V(\h)$ such that $|A|=|B|=n/2$ and $|E_*(A, B)\setminus E(\h)|\le \e n^k$, where $*$ denotes \emph{even} (respectively,  \emph{odd}) if $\h$ is $\e$-close to $\overline{\B}_{n,k}(A, B)$ (respectively, ${\B}_{n,k}(A, B)$).
For simplicity, we write $\B^*$ for $\overline\B_{n,k}(A,B)$ (respectively, for ${\B}_{n,k}(A, B)$) if $\h$ is $\e$-close to $\overline{\B}_{n,k}(A, B)$ (respectively, to ${\B}_{n,k}(A, B)$).

We observe that there are at most $\e_0 n$ vertices in $\h$ that are not $\e_1$-good with respect to $\B^*(A, B)$.
Indeed, recall that a vertex $v\in V(\h)$ is $\e_1$-good with respect to $\B^*(A, B)$ if $\deg_{\B^*(A, B)\setminus \h}(v)\le \e_1 n^{k-1}$. 
Since $|E_*(A, B)\setminus E(H)|\le \e n^k$, the number of  vertices which are not $\e_1$-good is at most $k\e n^k/(\e_1 n^{k-1}) = \e_0 n$.
\begin{lemma}\label{clm:A1B1}
There is a partition $A_1\cup B_1$ of $V(\h)$ with $|A_1|, |B_1|\ge (1/2-\e_0)n$ such that
 \begin{enumerate}[label=(\alph*)]
 \item every vertex $v$ is in at least $\frac15\binom{n-1}{k-1}$ majority edges (so there is no $\e_1$-bad vertex) with respect to $(A_1, B_1)$,\label{item:0}
 \item all but at most $(\e')^2 n^{k/2}$ $(k/2)$-sets are $\e'$-good with respect to $\B^*(A_1, B_1)$, and \label{item:I}
 \item at most $\e_0' n$ vertices are $\e_1$-medium with respect to $\B^*(A_1, B_1)$. \label{item:II}
 \end{enumerate}
\end{lemma}

\begin{proof}
Starting from the partition $A\cup B$, we obtain a new partition $A_1\cup B_1$ by moving all  vertices to the other part, that are not $\e_1$-good and  are contained in more minority edges than majority edges
 (that is, $ \deg_{\B^*(A, B)\cap \h}(v)\le \deg_{\h}(v) /2$). When a vertex $v$ is moved, all the edges of $\h$ that contain $v$ change parity. If we let $A'\cup B'$ denote the partition obtained from $A\cup B$ after moving $v$ to the other part, then 
$\deg_{\B^*(A', B')\cap \h}(v) = \deg_{\h}(v) - \deg_{\B^*(A, B)\cap \h}(v)$. Furthermore, since at most $\e_0 n$ vertices are not $\e_1$-good with respect to $\B^*(A, B)$, at most $\e_0 n$ are moved when deriving $A_1\cup B_1$. 
Therefore, for every vertex $v$ that is moved, we have
\begin{align*}
\deg_{\B^*(A_1, B_1)\cap \h}(v) &\ge \deg_{\h}(v) - \deg_{\B^*(A, B)\cap \h}(v) - (\e_0 n) n^{k-2} \\
&\ge \frac12 \left(\frac12-\e \right)\binom{n-1}{k-1} - \e_0 n^{k-1}\ge \frac15\binom{n-1}{k-1}.
\end{align*}
For every vertex $v$ that is not moved, we have
\[
\deg_{\B^*(A_1, B_1)\cap \h}(v) \ge \deg_{\B^*(A, B)\cap \h}(v) - (\e_0 n) n^{k-2} \ge \frac12 \left(\frac12-\e \right)\binom{n-1}{k-1} - \e_0 n^{k-1}\ge \frac15\binom{n-1}{k-1}.
\]
This proves ~\ref{item:0}. 

Since at most $\e_0 n$ vertices are moved, we have  $|A_1|, |B_1|\ge (1/2-\e_0)n$.
Moreover, we have  $|E_{*}(A_1, B_1)\setminus E(\h)|\le \e n^k+\e_0 n \binom{n-1}{k-1}\le \e_0 n^k$. 
This implies that the number of $(k/2)$-sets that are not  $\e'$-good is at most
\[
\binom{k}{k/2} \frac{\e_0 n^k}{\e' n^{k/2}} = \binom{k}{k/2} \frac{\e_0}{(\e')^3} (\e')^2 n^{k/2} \le (\e')^2 n^{k/2}
\]
as $\e$ is sufficiently small. This proves \ref{item:I}. By a similar calculation, we obtain that at most $\e_0' n$ vertices are $\e_1$-medium, proving \ref{item:II}. 
\end{proof}

Throughout the rest of the paper, whenever we use good, medium, bad, minority, majority etc., the underlying partition  referred to is always  $A_1\cup B_1$.
For brevity, we write $\overline{\B}:=\overline{\B}_{n,k}(A_1, B_1)$, $\B:={\B}_{n,k}(A_1, B_1)$ and $\B':=\B'_{n,k}(A_1, B_1)$.
Furthermore, let $N(S):=N_\h(S)$ denote the collection of $(k-|S|)$-sets $T$ in $V(\h)\setminus S$ such that $T\cup S\in E(\h)$.

\subsection{Lemmas for Step 1}\label{sec:step1}
In this section we collect all lemmas for Step 1.

\subsubsection{Cover all medium vertices}
The following lemma puts all medium vertices in a single $(k/2)$-path.

\begin{lemma}\label{clm:pm}
Let $M\neq\emptyset$ be the set of all $\e_1$-medium vertices in $\h$ and let $U$ be an arbitrary set of vertices of size at most $\e_0'n$.
Then there exists a $(k/2)$-path $\pp_M$ of length $4|M|-2$ in $\h\setminus U$ such that $\pp_M$ contains only edges in $E_*(A_1, B_1)$, $V(\pp)$ contains all $\e_1$-medium vertices and the ends of $\pp_M$ are $\e'$-good and both even or both odd depending on our choice.

\end{lemma}

\begin{proof}
We build the path with odd ends,
and the argument for even ends is the same. Fix a vertex $v$ and
let~$L_v$ be the hypergraph whose edges are  all $(k-1)$-sets $S$ which satisfy
\begin{itemize}
\item $S\cup \{v\}\in E_*(A_1, B_1)$, $S\cap U=\emptyset$,
\item $S$ contains no $\e_1$-medium vertex, and
\item all $(k/2)$-subsets of $S\cup \{v\}$ are $\e'$-good.
\end{itemize}
By Lemma~\ref{clm:A1B1}~\ref{item:II} and $|U|\leq \e_0' n$ there are at most $(\e_0' n+|U|) \binom{n-1}{k-2}\le \e_0' n^{k-1}$ $(k-1)$-sets
which fail to have the first two properties while by Lemma~\ref{clm:A1B1}~\ref{item:I} at most $(\e')^2 n^{k/2}\binom{n-k/2}{k/2-1}\le (\e')^2 n^{k-1}$ $(k-1)$-sets 
fail to have the third property.
As $v$ lies in at least $\frac15 \binom{n-1}{k-1}$ edges in $E_*(A_1, B_1)$ we infer therefore that $e(L_v)\ge \frac16  \binom{n-1}{k-1}$.

We greedily put vertices $v\in M$ in vertex-disjoint $(k/2)$-paths of length two, using majority edges and with $\e'$-good odd ends. 
This is possible since in each step there are at most $3(k/2)(4|M|-2)+|U|\le 10k \e_0' n$ chosen vertices and in turn, at least $\frac16  \binom{n-1}{k-1} - 10k \e_0' n \binom{n-1}{k-2}\ge \frac17  \binom{n-1}{k-1}$ edges in $L_v$ do not intersect these vertices. 
Among these edges we want to find two that share a $(k/2-1)$-set $T$ such that $T\cup \{v\}$ is odd if $*=\even$ and even otherwise.
This clearly yields the required length two path containing $v$. 
To the contrary, suppose that no such two edges exist and we shall derive a contradiction by counting $e(L_v)$.

Let $\T$ be the collection of $(k/2-1)$-sets $T$ such that $T\cup \{v\}$ is odd if $*=\even$ and $T\cup \{v\}$ is even otherwise.
For any $T\in \T$, let $\F_T$ be the family of edges of $L_v$ that contain $T$ and that do not intersect the chosen vertices. 
By our assumption, $\F_T$ must be intersecting, and thus, by the Erd\H{o}s--Ko--Rado theorem~\cite{EKR}, it has size 
at most $\binom{n-k/2-1}{k/2-1}$. 
Since there are at most $\binom {n-1}{k/2-1}$ choices for such~$T$, we obtain $\sum_{T\in \T}|\F_T|\le \binom {n-1}{k/2-1}\binom{n-k/2-1}{k/2-1}$.
Note that for any $e\in L_v$, if $e$ intersects both $A_1$ and $B_1$, then there exists a $(k/2-1)$-set $T\subseteq e$ such that $T\cup \{v\}$ satisfies the prescribed parity, namely, $T\in \T$ and $e$ is counted.
Therefore, the only members of $L_v$ possibly not counted in $\sum_{T\in \T}|\F_T|$ are the ones completely in $A_1$ or $B_1$.
The number of these, however, is at most $2\binom{(1/2+\e_0)n}{k-1}$, and thus the number of edges in $L_v$ not intersecting the chosen vertices is at most
\[
\binom {n-1}{k/2-1} \binom{n-k/2-1}{k/2-1} + 2\binom{(\frac12+\e_0)n}{k-1} 
< \frac17 \binom{n-1}{k-1},
\]
as $k\ge 6$ and $n$ is large. This is a contradiction.

It remains to connect these short paths to a single path $\pp_M$. This can be done by iteratively connecting two ends from two distinct paths by a $(k/2)$-set. This is possible since all the ends we have are $\e'$-good and the resulting path is not long. 
\end{proof}

\subsubsection{Building a bridge}

A crucial component of Step 1 is the construction of the \emph{bridges},  paths of constant size that overcome the parity issues. 
This is the only place where we use the minimum $d$-degree condition~\eqref{eq:md} in the proof.
The main tools are Lemmas~\ref{lem:deg_inter}--\ref{thm:deg_bridge}.
Recall that when we use good, medium, bad, minority, majority etc., the underlying partition  referred to is always  $A_1\cup B_1$.

The construction of the bridges depends on the parity of $2n/k$ and whether $\h$ is close to $\B$ or $\overline{\B}$. Note that when $n\in 2\mathbb N\setminus 4\mathbb N$
the hypergraph $\overline{\B}$ is not in $\h_{\ext}(n,k)$ so we do not need a bridge in this case.
Lemmas~\ref{clm:bridge},~\ref{clm:bridge2} and~\ref{clm:bridge3} handle each of the remaining three cases separately.

\begin{lemma}[Bridge for $\overline\B$ ($n\in k\mathbb N$)]\label{clm:bridge}
Suppose $n\in k\mathbb N$ and $\h$ is $\e$-close to $\overline{\B}$ and  satisfies~\eqref{eq:md}.
Assume that $|A_1|$ is odd. Then there exists a $(k/2)$-path $\pp_{\mathrm b}$ in $\h$ with $\e'$-good ends which has one of the following forms: $101$, $010$, $00100$, $11011$, or $001111100$.
\end{lemma}

\begin{proof}
Because $|A_1|$ is odd, we have $\overline{\B}\in \h_{\ext}(n,k)$ and odd edges are minority edges.
Suppose there is an  $\e'$-bad $(k/2)$-set $R$. Then $\deg_{\B\cap \h}(R) \ge \deg_{\h}(R) - \e'n^{k/2} \ge \frac14 \binom n {k/2}$, thus by Lemma~\ref{clm:A1B1}~\ref{item:I} the set 
$R$ forms an edge in $\h\cap \B$ with 
at least $\frac14 \binom n {k/2} - (\e')^2 n^{k/2} > \binom{n}{k/2-1}$ sets which are $\e'$-good. 
By the Erd\H{o}s--Ko--Rado theorem we can therefore find two disjoint $\e'$-good $(k/2)$-sets $S$ and $T$ among these neighbors of $R$ 
and  the $(k/2)$-path $SRT$ has the form 101  or 010,  as claimed.
We may  therefore assume that there is no $\e'$-bad $(k/2)$-set. 
By Lemma \ref{thm:deg_bridge}, $\h$ contains two odd edges that are either disjoint or sharing exactly $k/2$ vertices. 

If the former case occurs, then we partition these two disjoint 
edges arbitrarily into $R_1 S_1, R_2 S_2$ such that $R_1, R_2$ are odd $(k/2)$-sets (which thus are not $\e'$-bad). 
Further, by Lemma~\ref{clm:A1B1}~\ref{item:I} the number of $(k/2)$-sets that are $\e'$-medium is at most $(\e')^2 n^{k/2}$. Thus, 
we can find odd $\e'$-good $(k/2$)-sets $T_i\in N(R_i)$ for $i=1,2$ and even $\e'$-good $(k/2$)-sets $W_i\in N(S_i)$ for $i=1,2$ such that all chosen $(k/2)$-sets are disjoint. 
At last, since $T_1, T_2$ are $\e'$-good, we can pick an odd $(k/2)$-set $T\in N(T_1)\cap N(T_2)$ disjoint from all chosen sets. Now we get the path
\[
W_1~ S_1~ R_1~ T_1~ T~ T_2~ R_2~ S_2~ W_2
\]
which has the form 001111100 and $\e'$-good ends. 

In the latter case, when $\h$ contains two odd edges sharing exactly $k/2$ vertices, the path $S R T$ has the form 101 or 010. Since $S$ and~$T$ are not $\e'$-bad and 
the number of $(k/2)$-sets that at are $\e'$-medium is at most $(\e')^2 n^{k/2}$, we can find an $\e'$-good $(k/2)$-sets $S'\in N(S)$ (respectively, $T'\in N(T)$) with the same parity as $S$  (respectively, as $T$).
This yields a bridge $S' S R T T'$ of the form 00100 or 11011 with $\e'$-good ends. 
\end{proof}

\begin{lemma}[Bridge for $\B$ ($n\in k\mathbb{N}$)]\label{clm:bridge2}
Suppose $n\in k\mathbb N$ and $\h$ is $\e$-close to ${\B}$ and satisfies~\eqref{eq:md}.
Assume that $n/k-|A_1|$ is odd. Then there exists a $(k/2)$-path $\pp_{\mathrm b}$ in $\h$ with $\e'$-good ends, which has one of the following forms: $000$, $111$, $01110$, $10001$, or $100101001$.
\end{lemma}

\begin{proof}
The proof is very similar to the previous one.
Because $n/k-|A_1|$ is odd, ${\B}\in \h_{\ext}(n,k)$ and even edges are minority edges. 
Suppose there is an $\e'$-bad $(k/2)$-set $R$. Then $\deg_{\overline\B\cap \h}(R) \ge \deg_{\h}(R) - \e'n^{k/2} \ge \frac14 \binom n {k/2}$, 
thus by Lemma~\ref{clm:A1B1}~\ref{item:I} the set 
$R$ forms an edge in $\h\cap \overline\B$ with 
at least $\frac14 \binom n {k/2} - (\e')^2 n^{k/2} > \binom{n}{k/2-1}$ sets which are $\e'$-good. 
By the Erd\H{o}s-Ko-Rado theorem we  can therefore find two disjoint $\e'$-good $(k/2)$-sets $S$ and $T$ among these neighbors of $R$ 
and  the $(k/2)$-path $SRT$ has the form 000  or 111, as claimed.
We may therefore assume that there is no  $\e'$-bad $(k/2)$-set. By Lemma \ref{thm:deg_bridge}, $\h$ contains two even edges that are disjoint or  sharing $k/2$ vertices. 

In the former case, we partition these two edges arbitrarily into $R_1 S_1, R_2 S_2$ such that all of them are even $(k/2)$-sets (and thus none of them is $\e'$-bad). 
Further, by Lemma~\ref{clm:A1B1}~\ref{item:I} the number of $(k/2)$-sets that  are $\e'$-medium is at most $(\e')^2 n^{k/2}$. Thus  we can find odd $\e'$-good $(k/2$)-sets $T_i\in N(R_i)$ and $W_i\in N(S_i)$ for $i=1,2$ such that all chosen $(k/2)$-sets are disjoint. Lastly, since $T_1, T_2$ are $\e'$-good, we can pick an even $(k/2)$-set $T\in N(T_1)\cap N(T_2)$ disjoint from all chosen sets. Now we get a bridge
\[
W_1~ S_1~ R_1~ T_1~ T~ T_2~ R_2~ S_2~ W_2
\]
of the form 100101001 with $\e'$-good ends. 

In the latter case,  when $\h$ contains two even edges sharing exactly $k/2$ vertices, the path $S R T$ has the form  000 or 111. 
Since $S$ and~$T$ are not $\e'$-bad and 
the number of $(k/2)$-sets that are $\e'$-medium is at most $(\e')^2 n^{k/2}$, we can find $\e'$-good $(k/2)$-sets 
$S'\in N(S)$ (respectively, $T'\in N(T)$) which have  parity opposite to~$S$ (respectively, to $T$). 
So we get a bridge $S' S R T T'$ of the form 10001 or 01110 with $\e'$-good ends. 
\end{proof}

At last, consider the case $n\in \frac k2\mathbb N\setminus k \mathbb N$ and $\h$ is $\e$-close to ${\B}$.

\begin{lemma}[Bridge for $\B$ ($n\in \frac k2\mathbb N\setminus k \mathbb N$)]\label{clm:bridge3}
Suppose $n\in \frac k2\mathbb N\setminus k \mathbb N$ and $\h$ is $\e$-close to ${\B}$ and satisfies~\eqref{eq:md}.
There exists a $(k/2)$-path $\pp_{\mathrm b}$ in $\h$ with $\e'$-good ends which contains one or three even edges and satisfies the following.
\begin{itemize}
\item If $\lfloor n/k \rfloor - |A_1|$ is even, then the bridge has the form $1001$, $0110101110$ or $01101011010110$;
\item If $\lfloor n/k \rfloor - |A_1|$ is odd, then the bridge has the form $0110$, $1001010001$ or $10010100101001$.
\end{itemize}
\end{lemma}

The proof is not short so we give an outline here.
In the simplest case we find in $\h$ an even edge of the form 00 in case $\lfloor n/k \rfloor - |A|$ is even (i.e., the edge splits into two even $(k/2)$-sets)
or of the form 11 in case $\lfloor n/k \rfloor - |A|$ is odd.
 The Claim~\ref{clm:reduce} from below guarantees that one can then extend some even edge of this form from each of its end
  by an $\e'$-good $(k/2)$-set  and hence 
obtain a bridge of the form $1001$ or~$0110$. 

When such even edges do not exist we obtain by Claim~\ref{clm:reduce} and Fact~\ref{fact:evenk} from below a strong control over $\overline\B\cap \h$, the even edges of $\h$.
This results in Cases 
\ref{it:A}-\ref{it:D} in the proof
which quickly lead to a contradiction unless $\B'\in \h_{\ext}(n,k)$ and $\overline\B\cap \h=\h[A_1]\cup \h[B_1]$. 
This case requires more work and here we make  use of  Lemmas~\ref{lem:deg_inter} and~\ref{prop:e123} to find three edges $e_1, e_2, e_3\in \h[A_1]\cup \h[B_1]$ such that $e_1\cap (e_2\cup e_3)=\emptyset$ and $|e_2\cap e_3|\in \{0, k/2\}$. By suitably extending these configurations we  then obtain the desired bridge.

We start with the following auxiliary results. 
\begin{claim}\label{clm:reduce}
Suppose there is an even edge $e=S\cup T$ such that both $S$ and $T$ are even (respectively, odd) $(k/2)$-sets, then there is an even edge $e'$ (not necessarily distinct from $e$) which can be partitioned into two even (respectively, odd) $(k/2)$-sets that are not $\e'$-bad.
\end{claim}

\begin{proof}
We only prove the case when $S$ and $T$ are both even because the proof when they are both odd is identical.
For $1< i\le k/2$, let $\e_{i}:=k^{1/2} \e^{i/(2k)}$ and note that 
$\e' = \e_{k/2} \le\e_{k/2-1} \le\dots \le \e_1$. 
Suppose that~$S$ or $T$ is $\e'$-bad.
Let $\ell$ be the minimum integer such that there exists an $\epsilon_\ell$-bad $\ell$-set $L$. 
By Lemma~\ref{clm:A1B1}~\ref{item:0} there is no $\e_1$-bad vertex, thus $\ell>1$.
Further, since $S$ or $T$ is $\e_{k/2}$-bad, we have $\ell\leq k/2$.
Let $\ell_1=\lfloor \ell/2\rfloor$ and $\ell_2=\lceil \ell /2 \rceil$. 
We split $L$ arbitrarily into $L_1$ and $L_2$ of order $\ell_1$ and $\ell_2$, respectively.

Since $L$ is $\e_{\ell}$-bad, we have $\deg_{\B\cap \h}(L)\le \e_\ell n^{k-\ell}$. 
Since $|A_1|, |B_1|\ge (1/2-\e_0)n$, we have $\delta_{k-1}(\B)\ge (1/2-\e_0)n - k$, which, by Proposition~\ref{prop:seq}, implies that $\delta_\ell(\B)\ge (1/2-2\e_0) \binom{n-\ell}{k-\ell}$.
By the minimum degree condition of $\h$ and Proposition~\ref{prop:seq}, we have $\deg_{\h}(L)\ge (1/2-2\e_0) \binom{n-\ell}{k-\ell}$ and thus $\deg_{\overline\h}(L)\le (1/2+2\e_0) \binom{n-\ell}{k-\ell}$. Consequently, by $\deg_{\B\cap \h}(L)\le \e_\ell n^{k-\ell}$, we have 
\begin{equation}
\label{eq:315}
\deg_{\overline\B\cap \overline\h}(L) = \deg_{\overline\h}(L) - \deg_{\B\cap \overline\h}(L) \le (1/2+2\e_0) \binom{n-\ell}{k-\ell} - (\deg_{\B}(L) - \e_\ell n^{k-\ell}) < 2\e_\ell n^{k-\ell},
\end{equation}
where  the last inequality follows from the choice of $\e_\ell$.

Let $m(L_i)$, $i\in[2]$, denote the number of $(k/2-\ell_i)$-sets $L_i'$ such that $L_1\cup L_i'$ is an even $(k/2)$-set and not $\e'$-bad.
We  then claim that $m(L_i)\leq\e_{\ell_i} n^{k/2-\ell_i}$ for some $i\in[2]$.
Otherwise, the number of $(k-\ell)$-sets $L_1'\cup L_2'$ such that both $L_1\cup L_1'$ and $L_2\cup L_2'$ are even $(k/2)$-sets  and not $\e'$-bad is at least
$m(L_1)m(L_2)>(\e_{\ell_1} \e_{\ell_2} - o(1)) n^{k-\ell}=(k^{1/2}\e_{\ell} - o(1)) n^{k-\ell} >2\e_\ell {n}^{k-\ell}$. By \eqref{eq:315}, one of these $(k-\ell)$-sets lies in $N(L)$, 
that is, $e'=(L_1\cup L_1')\cup (L_2\cup L_2')$ is the desired even edge and we are done.

Without loss of generality, assume that $m(L_1) \leq \e_{\ell_1} n^{k/2-\ell_1}$.
Then
\[
\deg_{\B\cap \h}(L_1)\le \binom{n}{k/2-\ell_1} \e' n^{k/2} + \e_{\ell_1} n^{k/2-\ell_1} \binom{n}{k/2}< \e_{\ell_1} n^{k-\ell_1},
\]
which means that $L_1$ is $\e_{\ell_1}$-bad, contradicting the minimality assumption on $\ell$.
\end{proof}

\begin{fact}\label{fact:evenk}
Fix an even $k$-set $X\subset A_1\cup B_1$ with $x$ vertices in $A_1$ (thus $x$ is even).
\begin{enumerate}[label=$(\roman*)$]
\item If $0< x < k$, then $X$ can be partitioned into two odd $(k/2)$-sets. \label{item:evenk1}
\item If $0\le x < k$, then $X$ can be partitioned into two even $(k/2)$-sets. \label{item:evenk2}
\item If $x=k$, then $X$ can be partitioned into two   even $(k/2)$-sets when $k\in 4\mathbb{N}$ and into two
odd $(k/2)$-sets when $k\in 2\mathbb{N}\setminus 4\mathbb{N}$.  \qedhere  \label{item:evenk3} \qed
\end{enumerate}
\end{fact}



Now we are ready to prove Lemma~\ref{clm:bridge3}.

\begin{proof}[Proof of Lemma~\ref{clm:bridge3}]

Note that \eqref{eq:md} implies that $\delta_d(\h)> \delta_d(\B')$ if $\B'\in \h_{\ext}(n,k)$ and $\delta_d(\h) > \delta_d(\B)$ otherwise. Further,
by possibly swapping $A_1$ and $B_1$ we may assume that the partition $(A_1, B_1)$ satisfies $\deg_{\B}(S)\le \deg_{\B}(T)$ for any $S\in \binom{A_1}{d}$ and any $T\in \binom{B_1}{d}$.
We note that this implies 
\begin{enumerate}[label=$(\alph*)$]
\item if $\overline\B\cap \h =\h[A_1]\cup \h[B_1]$, then $\h[A_1]\neq\emptyset$. \label{item:A1}
\end{enumerate}
Indeed, otherwise $\overline\B\cap \h=\h[B_1]$ and
since $\deg_{\B}(S)\le \deg_{\B}(T)$ for any $S\in \binom{A_1}{d}$ and any $T\in \binom{B_1}{d}$, we can choose a $d$-set  $S\not\subseteq B_1$ such that $\deg_\B(S)=\delta_d(\B)$.
Since 
\[
\deg_\h(S)\ge \delta_d(\h)> \delta_d(\B)=\deg_\B(S),
\] 
the set $S$ must be contained in an even edge that intersects $A_1$, contradicting $\overline{\B} \cap \h = \h[B_1]$. \medskip

Suppose now that $\lfloor n/k \rfloor - |A_1|$ is even and  that $\h$ contains an even edge that can be split into two even $(k/2)$-sets. 
Then, by Claim~\ref{clm:reduce}, there is an even edge that can be split into two even $(k/2)$-sets $L_1, L_2$ that are not $\e'$-bad.
We can then pick $\e'$-good odd $(k/2)$-sets $L_3, L_4$ such that $L_3 L_1 L_2 L_4$ is a $(k/2)$-path of the form 1001 with good ends and we are done. 
In the case when $\lfloor n/k \rfloor - |A_1|$ is odd and $\h$ contains an even edge that can be split into two odd $(k/2)$-sets
an analogous argument applies and yields a $(k/2)$-path of the form 0110 with good ends.

Thus we may assume that such even edges do not exist and by Claim~\ref{clm:reduce}  and Fact~\ref{fact:evenk} conclude that $\h$ satisfies the following.
\begin{enumerate}[label=(\Alph*)]
\item\label{it:A} $\overline\B\cap \h=\emptyset$ when $k\in 4\mathbb{N}$ and $\lfloor n/k \rfloor - |A_1|$ is even;
\item\label{it:B} $\overline\B\cap\h =\h[A_1]$ when $k\in 2\mathbb{N}\setminus 4\mathbb{N}$ and $\lfloor n/k \rfloor - |A_1|$ is even;
\item\label{it:C} $\overline\B\cap\h =\h[A_1]\cup \h[B_1]$ when $k\in 4\mathbb{N}$ and $\lfloor n/k \rfloor - |A_1|$ is odd;
\item\label{it:D} $\overline\B\cap\h =\h[B_1]$ when $k\in 2\mathbb{N}\setminus 4\mathbb{N}$ and $\lfloor n/k \rfloor - |A_1|$ is odd.
\end{enumerate}

The Cases \ref{it:A} and \ref{it:D} immediately contradict~\ref{item:A1} and in the following we deal with the 
Cases \ref{it:B} and~\ref{it:C} for which we have $\B'\in \h_{\ext}(n,k)$.
Thus, it suffices to consider the case $\B'\in \h_{\ext}(n,k)$
together with the assumption $\overline\B\cap \h =\h[A_1]\cup \h[B_1]$, which we do in the following.

\medskip

We first claim that none of $L\in \binom{A_1}{k/2}\cup \binom{B_1}{k/2}$ is $\e'$-bad.
Fix $L\in \binom{A_1}{k/2}\cup \binom{B_1}{k/2}$.
Since there is no even $(i,{k-i})$-edge, $0< i<k$ (recall that an edge $e$ is an $(i,j)$-edge if $|e\cap A_1|=i$), which contains $L$, we know that
\begin{equation}\label{eq:degA1}
\deg_{\overline\B\cap \h}(L) \le \binom{|A_1|}{k/2}\le (1/2+\e_0)^{k/2}\binom{n}{k/2}\le \frac13\binom{n}{k/2}
\end{equation}
as $k\ge 4$ and $\e_0$ is small enough.
Together with $\delta_{k/2}(\h)\ge (1/2 - \e) \binom{n}{k/2}$, we infer that $\deg_{\B\cap \h}(L)=\deg_\h(L) - \deg_{\overline\B\cap \h}(L)\ge \frac17\binom{n}{k/2}$, i.e.,~$L$ is not $\e'$-bad.

Next we show that \emph{there exist $e_1, e_2, e_3\in \h[A_1]\cup \h[B_1]$ such that $e_1\cap (e_2\cup e_3)=\emptyset$ and $|e_2\cap e_3|\in \{0, k/2\}$}.
Let $S_0$ be a $d$-set such that $\deg_{\B'}(S_0)=\delta_d(\B')$.
If $S_0$ intersects both $A_1$ and $B_1$, then as $\deg_\h(S_0)> \delta_d(\B')=\deg_{\B'}(S_0)$, the set $S_0$ must be contained in an even edge that intersects both $A_1$ and~$B_1$, contradicting $\overline\B\cap \h =\h[A_1]\cup \h[B_1]$.
Thus, $S_0\subseteq A_1$ or $S_0\subseteq B_1$. 

First, assume $S_0 \subseteq A_1$. 
Since $\deg_{\B'}(S_0)=\delta_d(\B')$, every $d$-set $S \subseteq A_1$ satisfies
\[
\deg_\h(S) > \delta_d(\B') = \deg_{\B'}(S_0)\ge \deg_\B(S_0)+\delta_d(\mathcal S_{|A_1|, k}) = \deg_\B(S)+\delta_d(\mathcal S_{|A_1|, k}).
\]
Since $\overline\B\cap \h =\h[A_1]\cup \h[B_1]$, any even edge containing $S$ must be entirely in $A_1$.
Consequently, $\deg_{\h[A_1]}(S)\ge \deg_\h(S) - \deg_\B(S) > \delta_d(\mathcal S_{|A_1|, k})$, implying that $\delta_d(\h[A_1])> \delta_d(\mathcal S_{|A_1|, k})$.
Applying Lemma~\ref{prop:e123} to $\h[A_1]$ gives the desired $e_1, e_2, e_3$.

Second, assume $S_0\subseteq B_1$. 
In this case $\delta_d(\B')$ must be attained by every $d$-set $S\subseteq B_1$. 
It follows that for any such $S$, $\deg_\h(S)\ge \delta_d(\B')+1= \deg_{\B'}(S)+1$, and consequently, $\deg_{\h\cap \overline\B}(S)\ge 1$.
Since all even edges containing $S$ must be entirely in $B_1$, we have
\[
e(\h[B_1])\ge \binom{|B_1|}{d} / \binom{d}{k/2} = \Omega(n^{d})= \Omega(n^{k/2}).
\]
Thus, we can find two edges $e_2, e_3\subseteq B_1$ such that $|e_2\cap e_3|\in \{0, k/2\}$ by applying Lemma \ref{lem:deg_inter} to $\h[B_1]$. 
By~\ref{item:A1}, there is an edge $e_1\in \h[A_1]$.
Since $e_1\cap (e_2\cup e_3)=\emptyset$, we obtain the desired $e_1$, $e_2$ and $e_3$. 
\medskip

Finally, we construct the bridge from $e_1$, $e_2$ and $e_3$.
We will only show the case when $k\in 2\mathbb{N}\setminus 4\mathbb{N}$ and $\lfloor n/k \rfloor - |A_1|$ is even because the case when $k\in 4\mathbb{N}$ and $\lfloor n/k \rfloor - |A_1|$ is odd is identical after exchanging even with odd and exchanging 0 with 1.
If $|e_2\cap e_3|=0$ then we split $e_i$,  $i\in [3]$, into two (odd) $(k/2)$-sets $L_i, L_i'$. 
Since $L_i, L_i'$ are not $\e'$-bad, we can find disjoint $\e'$-good even $(k/2)$-sets $S_i\in N(L_i), S_i'\in N(L_i)$ for $i\in [3]$. 
Finally, we pick two odd $(k/2)$-sets $T_1\in N(S_1')\cap N(S_2)$ and $T_2\in N(S_2')\cap N(S_3)$ such that they are all disjoint and disjoint from all chosen sets. This yields the path
\[
S_1~L_1~L_1'~S_1'~T_1~S_2~L_2~L_2'~S_2'~T_2~S_3~L_3~L_3'~S_3'
\]
which has  the form 01101011010110. Otherwise let $e_2\cap e_3=L_0$ and  we split $e_1$ into $L_1\cup  L_2$ 
and let $L_3=e_2\setminus e_3$ and $L_4=e_3\setminus e_2$. Similarly, we pick $\e'$-good even $(k/2)$-sets 
$S_i\in N(L_i)$ for $i\in [4]$ and an odd $(k/2)$-set $T\in N(S_2)\cap N(S_3)$ such that they are all disjoint and disjoint from all chosen sets. This yields the path
\[
S_1~L_1~L_2~S_2~T~S_3~L_3~L_0~L_4~S_4
\]
which has the form 0110101110.
\end{proof}

\subsection{Proof of Theorem \ref{lemE}}
We now prove  Theorem \ref{lemE}, following the overview from Section~\ref{sec:overview}. 
In particular, we use the partition $A_1\cup B_1$ obtained in Lemma~\ref{clm:A1B1} and
as mentioned in the overview, the proof of  Theorem \ref{lemE} splits into several cases depending 
on whether $n\in \frac k2\mathbb N\setminus k \mathbb N$ or not and whether 
$\h$ is $\e$-close to $\overline{\B}$ or ${\B}$.
Each case may be further split
depending on the parity of the ends of $\pp$, the path to be established in Step 1. 
We recall also that the binary representations of edges from $\overline{\B}$ in a $(k/2)$-path are~$00$ or $11$ while those of edges in~$\B$ are $01$ or $10$.

\subsubsection{The case $n\in k\mathbb N$ and $\h$ is $\e$-close to $\overline{\B}$}
\label{sec:case1}

In this case $n$ is even since $k$ is even. 

\medskip
\noindent\textbf{Step 1.} 
We first build a short path $\pp$ that contains the bridge from Lemma~\ref{clm:bridge} and all medium vertices.

\begin{claim}\label{clm:shortpath}
There exists a $(k/2)$-path $\mathcal P$ in $\h$ such that 
\begin{itemize}
\item $|V(\pp)|\le 3k\e_0' n$,
\item $V(\pp)$ contains all $\e_1$-medium vertices,
\item the ends of $\mathcal P$ are  $\e'$-good $(k/2)$-sets with the same parity,
\item $\pp$ has an odd length and $|A_1\setminus V(\pp)|$ is even.
\end{itemize}
\end{claim}

\begin{proof}
We separate cases based on the parity of $|A_1|$.

\medskip
\noindent\textit{Case 1.} $|A_1|$ is odd (and thus $\overline{\B}\in \h_{\ext}(n,k)$). 
\medskip

We apply Lemma \ref{clm:bridge} to obtain the path $\pp_{\mathrm b}$ of even length, whose ends are either both odd (i.e., $\pp$ has the binary representation $101$ or $11011$) or both even
 (i.e., $\pp$ is of the form  $010$, $00100$, or $001111100$). 
Further, let $M$ denote the set of $\e_1$-medium vertices and apply Lemma~\ref{clm:pm}  with $U=V(\pp_{\mathrm b})$ 
to obtain the path $\pp_M$, which  has even length $4|M|-2$, is disjoint from $U$ and covers all vertices in $M$. 
Moreover, Lemma~\ref{clm:pm} allows us to choose  the parity of the ends of~$\pp_M$ and we choose it to be the same as the ends of $\pp_{\mathrm b}$. In particular, if they are odd (respectively, even)
then the binary representation of~$\pp_M$ is an all 1 string of odd length due to its even length (respectively, a constant 0 string).

As all ends are $\e'$-good and have the same parity
we can connect one end of $\pp_{\mathrm b}$ with one end of $\pp_M$ by a $(k/2)$-set with the same parity as these ends.
We extend the path by one edge to obtain the path~$\pp$ odd length while keeping the ends  $\e'$-good and of the same parity. 
Note that $|V(\pp)|\le (9+4\e_0' n+1)k/2\le 3k\e_0' n$ and that there is an odd number of $1$'s in the binary representation of $\pp$. Thus
$|V(\pp)\cap A_1|$ is odd, 
$|A_1\setminus V(\pp)|$ is even and $\pp$ is the desired $(k/2)$-path. 

\bigskip
\noindent\textit{Case 2.} $|A_1|$ is even (and thus $\overline{\B}\notin \h_{\ext}(n,k)$). 
\medskip

We find the path $\pp_M$ with even ends by applying Lemma~\ref{clm:pm} with $U=\emptyset$.
We extend the path by one more edge to make its length odd and the resulting path $\pp$ satisfies all the requirements of Claim~\ref{clm:shortpath}. 
\end{proof}

\medskip
\noindent\textbf{Step 2.} Let $\pp$ be the $(k/2)$-path obtained from Step 1, Claim~\ref{clm:shortpath}, with ends denoted by  $L$ and $S$. Let $A_1'=A_1\setminus V(\pp)$ 
and  $B_1'=B_1\setminus V(\pp)$ and note that $|A_1'|$ is even by Claim~\ref{clm:shortpath} .  
We will extend $\pp$ to a Hamilton $(k/2)$-cycle by applying Lemma \ref{lem:finish}  or 
 Lemma \ref{lem:finish1} depending on $k$ and the parity of $L$ and $S$. Before being able to do so we need to make some adjustments to the partitions.

\medskip
\noindent\textit{Case $(i)$}. 
The sets $L, S$ are even and $k\in 4\mathbb N$. 

Let $b\equiv|A_1'|$ mod $k/2$ such that $0\le b<k/2$ and note that $b$ is even as  $|A_1'|$ and $k/2$ are both even. 
We  pick an $\e'$-good $(b, k/2-b)$-set $L'\in N(L)$ in  $A_1'\cup B_1'$ and thus have $|A_1'\setminus L'|\in \frac k2 \mathbb N$ and $|B_1'\setminus L'|\in \frac k2 \mathbb N$. 
Next we pick  $\e'$-good $(0, k/2)$-sets $L_1, L_2$ and $\e'$-good $(k/2, 0)$-sets $S_1, S_2$ from $(A_1'\cup B_1') \setminus L'$, all disjoint and such that $L' L_1, SS_1, L_2 S_2\in E(\h)$. 

Let us verify the assumptions of Lemma \ref{lem:finish}. 
For any $v\in A_1'\setminus  L'$, we have $\deg_{\overline\B\setminus \h}(v)\le \e_1 n^{k-1}$ since $v$ is $\e_1$-good. 
Further, as $n_1:=|A_1'\setminus  L'|=(1- o(1)) n/2$, we have
\[
\deg_{\overline{\h[A_1]}}(v)\le \deg_{\overline\B\setminus \h}(v)\le \e_1 n^{k-1}\le 2^k \e_1 n_1^{k-1},
\]
where $\overline{\h[A_1]}$ represents the complement of $\h[A_1]$ (on the vertex set $A_1$). 
Similarly, $\deg_{\overline{\h[A_1]}}(S_i)\le 2^k \e' n_1^{k/2}$ for $i=1,2$.
That is, every $v\in A_1$ and $S_1, S_2$ are $(2^k\e_1)$-good with respect to $\h[A_1'\setminus  L']$ and
we apply Lemma \ref{lem:finish} on $\h[A_1'\setminus L']$ with $\a_0=2^k\e_1$ and sets $S_1, S_2$ and obtain a Hamilton path $\pp_1$ with ends $S_1, S_2$. 
Similarly, we apply Lemma \ref{lem:finish} to $\h[B_1'\setminus  L']$ with $\a_0=2^k\e_1$ and sets $L_1, L_2$ and obtain a Hamilton path $\pp_2$ with ends $L_1, L_2$. 
This yields the Hamilton $(k/2)$-cycle 
\[S~\pp~L~L'~L_1~\pp_2~L_2~S_2~\pp_1~S_1~S.\]

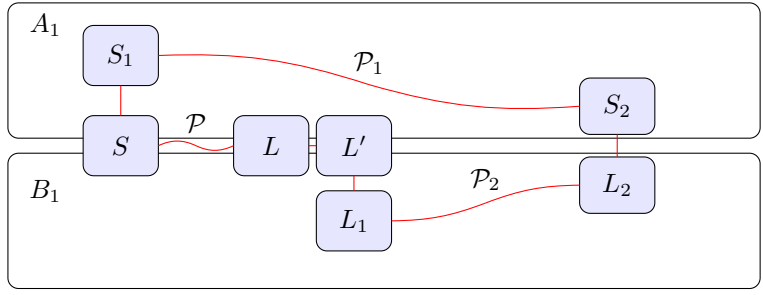
\begin{figure}[h]
\begin{center}
\begin{tikzpicture}
[inner sep=2pt,
   vertex/.style={circle, draw=blue!50, fill=blue!50},
   ]
\begin{pgfonlayer}{bg}    
\draw[rounded corners] (0,0) rectangle (10, 1.8);
\draw[rounded corners] (0,-2) rectangle (10, -0.2);
\end{pgfonlayer}
\draw[rounded corners, fill=blue!10] (1,0.7) rectangle (2,1.5); 
\draw[rounded corners, fill=blue!10] (1,-0.5) rectangle (2,0.3); 
\draw[rounded corners, fill=blue!10] (3,-0.5) rectangle (4,0.3); 
\draw[rounded corners, fill=blue!10] (4.1,-0.5) rectangle (5.1,0.3);  
\draw[rounded corners, fill=blue!10] (4.1,-1.5) rectangle (5.1,-0.7); 
\draw[rounded corners, fill=blue!10] (7.6,-1) rectangle (8.6,-0.25); 
\draw[rounded corners, fill=blue!10] (7.6,0.05) rectangle (8.6,0.8); 

\node at (0.5,1.5) {$A_1$};
\node at (0.5,-0.7) {$B_1$};
\node (S1) at (1.5,1.1) {$S_1$};
\node (S) at (1.5, -0.1) {$S$};
\node (L) at (3.5, -0.1) {$L$};
\node (LL) at (4.6, -0.1) {$L'$};
\node (L1) at (4.6, -1.1) {$L_1$};
\node (L2) at (8.1, -0.625) {$L_2$};
\node (S2) at (8.1, 0.425) {$S_2$};
\draw[color=red]  (1.5, 0.3) -- (1.5, 0.7); 
\draw[color=red]  (4, -0.1) -- (4.1, -0.1); 
\draw[color=red]  (4.6, -0.5) -- (4.6, -0.7); 
\draw[color=red]  (8.1, 0.05) -- (8.1, -0.25); 
\draw[color=red] (2,-0.1) to [bend left=25] (2.5, -0.1);
\draw[color=red] (2.5,-0.1) to [bend right=25] (3, -0.1);
\node at (2.5, 0.2) {$\pp$};
\draw[color=red] (2,1.1) to [bend left=10] (4.8, 0.7125); 
\draw[color=red] (4.8,0.7125) to [bend right=10] (7.6, 0.425);
\node at (4.8, 1) {$\pp_1$};
\draw[color=red] (5.1,-1.1) to [bend right=10] (6.35, -0.8625);
\draw[color=red] (6.35,-0.8625) to [bend left=10] (7.6, -0.625);
\node at (6.35, -0.55) {$\pp_2$};
\end{tikzpicture}

\caption{An illustration of the Hamilton cycle in Case $(i)$}
\end{center}
\end{figure}

\medskip
\noindent\textit{Case $(ii)$}. The sets  $L, S$ are even and $k\in 2\mathbb N\setminus 4\mathbb N$. 

Let $b\equiv|A_1'|$ mod $(k/2-1)$ such that $0\le b<k/2-1$. 
Then $b$ is even  as  $|A_1'|$ and  $k/2 -1$  are both even.  
We  pick an $\e'$-good $(b, k/2-b)$-set $L'$  from $A_1'\cup B_1'$ such that $LL'\in E(\h)$ and consequently $|A_1'\setminus L'|\in (\frac k2-1) \mathbb N$. 
Next we pick $\e'$-good $(0, k/2)$-sets $L_1, L_2$ and $\e'$-good $(k/2-1, 1)$-sets $S_1, S_2$ from $(A_1'\cup B_1') \setminus L'$, all disjoint and such that $L' L_1, SS_1, L_2 S_2\in E(\h)$. 
Let $Y$ be an arbitrary subset of $B_1'\setminus  (L'\cup L_1\cup L_2)$ of order $\frac2{k-2}|A_1'\setminus L'|$ which contains $S_1\cap B_1$ and $S_2\cap B_1$.
Since
\[
|(A_1'\setminus L')\cup Y| = \frac {k-2}2 |Y| +|Y| = \frac k2 |Y| \in \frac k2 \mathbb N,
\]
we infer that $|B_1'\setminus (L'\cup Y)|\in \frac k2 \mathbb N$.
We apply Lemma \ref{lem:finish1} on $\h[(A_1'\setminus L')\cup Y]$ with $\a_0=2^k\e_1$ and sets $S_1, S_2$ and obtain a Hamilton path $\pp_1$ with ends $S_1, S_2$. 
Note that $|B_1'\setminus (L'\cup Y)|\geq (1-\frac{2}{k-2}-o(1))\frac n2\ge \frac n4$ and
it is easy to check that every vertex is $\a_0$-good, with $\a_0=\sqrt \e_1$,  and that $L_1$ and $L_2$ are both $\a_0$-good  with respect to $\K^k(B_1'\setminus (L'\cup Y))$.
We apply Lemma \ref{lem:finish} on $\h[B_1'\setminus (L'\cup Y)]$ with $\a_0=\sqrt\e_1$ and sets $L_1, L_2$ to obtain a Hamilton path $\pp_2$ with ends $L_1, L_2$. Thus, we get a Hamilton $(k/2)$-cycle
\[S~\pp~L~L'~L_1~\pp_2~L_2~S_2~\pp_1~S_1~S.\]

\medskip
\noindent\textit{Case $(iii)$}. The sets $L, S$ are odd and $k\in 2\mathbb N\setminus 4\mathbb N$. 

This case becomes Case $(ii)$ after we exchange $A_1$ and $B_1$ (thus $L$ and $S$ become even).

\medskip
\noindent\textit{Case $(iv)$.} The sets $L, S$ are odd and $k\in 4\mathbb N$. 

Without loss of generality, assume $|B_1'| \ge |A_1'|$ and let
\[
b:=|B_1'| - |A_1'|=n-|V(\pp)| - 2|A_1'|.
\]
As  $|A_1'|$ is even, $\pp$ has odd length and $n\in k\mathbb N$  we have $b\in 4\mathbb N$. 
Further, as $|A_1|, |B_1|\geq (1/2-\e_0) n$ and $|V(\pp)|<3k\e_0'n$, we have
\[
b\le ||B_1|-|A_1|| + |V(\pp)| \le 2\e_0 n + 3k\e_0' n \le 4k \e_0' n.
\]
To balance out the sizes of $A_1'$ and $B_1'$ we  extend $\pp$ slightly, using more vertices from $B_1'$ while keeping the main properties of $\pp$.
If $k/4-1$ is odd, we greedily extend $\pp$ from $L$ by $b/2$ $\e'$-good $(k/4-1, k/4+1)$-sets and denote the path obtained by $\pp'$. 
Otherwise $k/4-1$ is even and we greedily extend the path $\pp$ from $L$ by $b/4$ $\e'$-good $(k/4-2, k/4+2)$-sets (note that $k\ge 12$ in this case). 
We denote the resulting path by $\pp'$ if $b/4$ is even. 
Otherwise, we extend the path by one more $\e'$-good $(k/4, k/4)$-set and let the resulting path be $\pp'$. Note that the above process is possible since all sets involved are odd and $\e'$-good. 
Let  $L'$ be the new end of $\pp'$ and let $A_1''=A_1\setminus V(\pp')$ and $B_1''=B_1\setminus V(\pp')$. 
Then $|A_1''| = |B_1''|=:m$ and $|V(\pp')|\le |V(\pp)| + (b/2)\cdot k/2+k/2\le 3k\e_0'n + k^2 \e_0' n+k/2\le 2k^2\e_0' n$. 
By definition, $\pp'$ has an odd length and consequently $|V\setminus V(\pp')|\in k\mathbb N$ and $m\in (k/2)\mathbb N$.

Next we pick $\e'$-good $(k/2-1, 1)$-sets $L_1, L_2$ and $\e'$-good $(1, k/2-1)$-sets $S_1, S_2$ from $V\setminus V(\pp')$ such that $L' L_1, SS_1, L_2 S_2\in E(\h)$. 
Let $X$ be an arbitrary subset of $A_1''\setminus (L_1\cup L_2)$ of order $2m/k$ containing $S_1\cap A_1$ and $S_2\cap A_1$, 
and let $Y$ be an arbitrary subset of $B_1''\setminus (S_1\cup S_2)$ of order $2m/k$ containing $L_1\cap B_1$ and $L_2\cap B_1$.
Then $|A_1''\setminus X|=|B_1''\setminus Y|\in (\frac k2-1) \mathbb N$ and
$|A_1''\setminus X|$, $|B_1''\setminus Y|= (1-\frac{2}{k}-o(1))\frac n2\ge \frac 38 n$. 
We apply Lemma~\ref{lem:finish1} on $\h[(B_1''\setminus Y)\cup X]$ with $\a_0=\sqrt\e_1$ and sets $S_1, S_2$ and obtain a Hamilton path $\pp_1$ with ends $S_1, S_2$
and apply Lemma \ref{lem:finish1} on $\h[(A_1''\setminus X)\cup Y]$ with $\a_0=\sqrt\e_1$ and sets $L_1, L_2$ and obtain a Hamilton path $\pp_2$ with ends $L_1, L_2$. 
This yields the Hamilton $(k/2)$-cycle 
\[
S~ \pp'~ L'~ L_1~ \pp_2~ L_2~ S_2~ \pp_1~ S_1~ S.
\]

\begin{figure}[h]
\begin{center}
\begin{tikzpicture}
[inner sep=2pt,
   vertex/.style={circle, draw=blue!50, fill=blue!50},
   ]
\begin{pgfonlayer}{bg}    
\draw[rounded corners] (-1,0) rectangle (9, 1.8);
\draw[rounded corners] (-1,-2) rectangle (9, -0.2);
\end{pgfonlayer}
\draw[rounded corners, fill=blue!10] (-0.2,-0.65) rectangle (0.8,0.15); 
\draw[rounded corners, fill=blue!10] (1,-0.5) rectangle (2,0.3); 
\draw[rounded corners, fill=blue!10] (2.8,-0.5) rectangle (3.8,0.3); 
\draw[rounded corners, fill=blue!10] (4.1,-0.35) rectangle (5.1,0.45); 
\draw[rounded corners, fill=blue!10] (5.9,-0.35) rectangle (6.9,0.45); 
\draw[rounded corners, fill=blue!10] (7.2,-0.65) rectangle (8.2,0.15); 
\draw (3.9,-0.2) to [bend right=60] (7.1, -0.2); 
\draw (-1,1) to [bend left=30] (0.9, 0); 
\draw (9,1) to [bend right=30] (7, 0); 

\node at (-0.5,1.5) {$A_1$};
\node at (-0.5,-0.7) {$B_1$};
\node (S1) at (0.3,-0.3) {$S_1$};
\node (S) at (1.5, -0.1) {$S$};
\node (L) at (3.3, -0.1) {$L'$};
\node at (5.5, -0.7) {$B_1'$};
\node at (-0.6, 0.5) {$A_1'$};
\node at (8.6, 0.5) {$A_1'$};
\node (L1) at (4.6, 0.05) {$L_1$};
\node (L2) at (6.4, 0.05) {$L_2$};
\node (S2) at (7.7, -0.3) {$S_2$};
\draw[color=red]  (0.8, -0.1) -- (1, -0.1); 
\draw[color=red]  (3.8, -0.1) -- (4.1, -0.1); 
\draw[color=red]  (6.9, -0.1) -- (7.2, -0.1); 
\draw[color=red] (2,-0.1) to [bend left=25] (2.4, -0.1);
\draw[color=red] (2.4,-0.1) to [bend right=25] (2.8, -0.1);
\node at (2.5, 0.2) {$\pp'$};
\draw[color=red] (0.3,-0.65) to [bend right=20] (7.7, -0.65); 
\node at (4, -1.15) {$\pp_1$};
\draw[color=red] (4.6,0.45) to [bend right=20] (2.7, 0.825);
\draw[color=red] (2.7,0.825) to [bend left=20] (0.8, 1.2);
\draw[color=red] (0.8, 1.35) to [bend left=10] (7.6, 1.35);
\draw[color=red] (0.8, 1.2) to [bend left=70] (0.8, 1.35);
\draw[color=red] (7.6, 1.35) to [bend left=70] (7.6, 1.2);
\draw[color=red] (7.6,1.2) to [bend left=20] (7.1, 0.825);
\draw[color=red] (7.1,0.825) to [bend right=20] (6.6, 0.45);
\node at (4.2, 1.3) {$\pp_2$};
\end{tikzpicture}

\caption{An illustration of the Hamilton cycle in Case $(iv)$}
\end{center}
\end{figure}
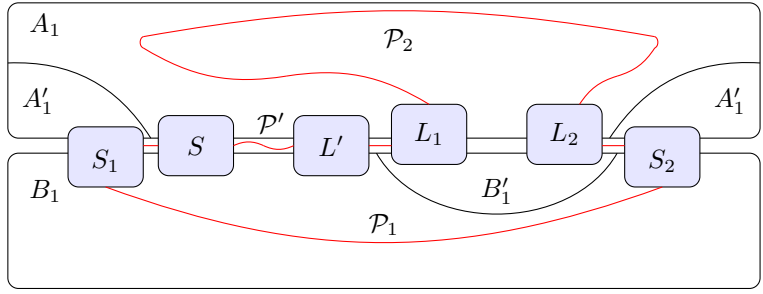

\subsubsection{The case $n\in \frac k2\mathbb N\setminus k \mathbb N$ and $\h$ is $\e$-close to ${\overline\B}$}

The proof of this subsection is almost identical to the one in Section~\ref{sec:case1}. 
As $\overline{\B}\notin \h_{\ext}(n,k)$, we do not need the bridge $\pp_{\mathrm b}$. 
We apply Lemma~\ref{clm:pm} to build a path that i) contains all medium vertices, ii) has even length, and iii) with good ends that have the same parity as $|A_1|$ (recall that Lemma~\ref{clm:pm} allows us to decide the parity of the ends).
Note that if we extend this path by one more good $(k/2)$-set, we obtain a path $\pp$ that satisfies all criteria of Claim~\ref{clm:shortpath}.
Step 2 is the same as in Section~\ref{sec:case1}.

\subsubsection{The case $n\in k \mathbb N$ and $\h$ is $\e$-close to ${\B}$} $\,$

\medskip
\noindent\textbf{Step 1.} 
We also find a short path $\pp$ as in the previous cases.
\begin{claim}\label{clm:shortpath2}
There exists a $(k/2)$-path $\mathcal P$ in $\h$ such that 
\begin{itemize}
\item $|V(\pp)|\le 3k\e_0' n$,
\item $V(\pp)$ contains all $\e_1$-medium vertices,
\item the ends of $\mathcal P$ are  $\e'$-good $(k/2)$-sets with different parities,
\item $\pp$ has an odd length and $|A_1\setminus V(\pp)|-\frac{n-|V(\pp)|}{k}$ is even.
\end{itemize}
\end{claim}

\begin{proof}
We separate cases based on the parity of $n/k-|A_1|$.

\medskip
\noindent\textit{Case 1.} $n/k-|A_1|$ is odd (and thus ${\B}\in \h_{\ext}(n,k)$).
\medskip

Let $\pp_{\mathrm b}$ be given by Lemma~\ref{clm:bridge2} which contains exactly two even edges. 
Let $M$ be the set of 
$\e_1$-medium vertices and we apply Lemma~\ref{clm:pm} with $U=V(\pp_{\mathrm b})$ to find a path $\pp_M$, which covers $M$, has  even length and 
such that its ends have the same parity as those of~$\pp_{\mathrm b}$. 
As the ends of~$\pp_M$ and $\pp_{\mathrm b}$  are $\e'$-good and have same parity we can connect these paths by one $(k/2)$-set whose parity is opposite to these ends.
Note that since both $\pp_M$ and $\pp_{\mathrm b}$ have even lengths, the length of the resulting path is also even. 

We extend the path by one more edge to make its length odd and denote the resulting path by $\pp$. Note that the ends of $\pp$ are $\e'$-good, have different parities 
and $|V(\pp)|\le (9+4\e_0' n+1)k/2\le 3k\e_0' n$. 
It remains to show that $|A_1\setminus V(\pp)| - \frac{n-|V(\pp)|}{k}$ is even. 
Let $t:=|V(\pp)|/k\in \mathbb N$.
Since all but two edges of~$\pp$ are odd,~$\pp$ contains $t+1$ (when  $\pp_{\mathrm b}$  has the form $111$ or $01110$) or $t-1$ 
(when  $\pp_{\mathrm b}$  has the form $000$, $10001$  or $100101001$) odd $(k/2)$-sets.
This implies that $|A_1\setminus V(\pp)|\equiv |A_1| - (t-1) \pmod 2$. Thus we have
\[
|A_1\setminus V(\pp)|\equiv |A_1| - (t-1) \equiv \frac nk -1- (t-1) = \frac nk - t = \frac{n-|V(\pp)|}{k} \pmod 2.
\]

\medskip
\noindent\textit{Case 2.} $n/k-|A_1|$ is even.
\medskip

Since ${\B}\notin \h_{\ext}(n,k)$, we do not need a bridge to correct the parity. 
We find the path $\pp_M$ with even ends by applying Lemma~\ref{clm:pm} with $U=\emptyset$, and then extend $\pp_M$ by one good odd $(k/2)$-set.
Denote the resulting path by $\pp$. So $\pp$ has an odd length and its ends have different parities. 
It remains to show that $|A_1\setminus V(\pp)|-\frac{n-|V(\pp)|}{k}$ is even. 
Note that by definition, $\pp$ contains $t:=|V(\pp)|/k$ odd $(k/2)$-sets, and thus
\[
|A_1\setminus V(\pp)|\equiv |A_1| -t \equiv \frac nk - t = \frac{n-|V(\pp)|}{k} \pmod 2. \qedhere
\]
\end{proof}

\medskip
\noindent\textbf{Step 2.}
Let $\pp$ be the $(k/2)$-path obtained from Step 1, Claim~\ref{clm:shortpath2}, with ends denoted by  $L$ and $S$ where~$L$ is odd and $S$ is even. Let $n_1=|A_1\setminus V(\pp)|$
and  $n_2=|B_1\setminus V(\pp)|$. Without loss of generality assume $n_2 \ge n_1$ and let $b:=n_2 - n_1$. 
We will extend $\pp$ to a Hamilton $(k/2)$-cycle by applying Lemma \ref{lem:finish2} and are thus required to make some adjustments before being able to do so.
 
Note that
\begin{equation}\label{eq:bbbb}
b=n_2 - n_1= n-|V(\pp)| - 2n_1= 2\left( \frac{n-|V(\pp)|}{k} - n_1 \right) + (k-2)\frac{n-|V(\pp)|}{k} .
\end{equation}
Moreover, by the definition of $A_1$ and $B_1$ and Claim~\ref{clm:shortpath2}, we have
\[
b\le ||B_1|-|A_1|| + |V(\pp)| \le 2\e_0 n + 3k\e_0' n \le 4k \e_0' n.
\]
We separate two cases according to the parity of $k/2$.

\medskip
\noindent\textit{Case $i)$}. $k\in 2\mathbb{N}\setminus 4\mathbb{N}$. 
By Claim \ref{clm:shortpath2}, $n_1 - \frac{n-|V(\pp)|}{k}$ is even. Note that $k-2\in 4\mathbb N$ and together with \eqref{eq:bbbb}, we have $b\in 4\mathbb N$.

Note that $k/4-1/2$ and $k/4-3/2$ are two consecutive integers and without loss of generality, assume that $k/4-1/2$ is even. 
We greedily extend the path $\pp$ from $L$ by $b/4$ $\e'$-good $(k/4-1/2, k/4+1/2)$-sets and $b/4$ $\e'$-good $(k/4-3/2, k/4+3/2)$-sets alternately. 
This process is possible since all edges involved are odd. 
Let the resulting path be $\pp'$ and denote its new end by $L'$.
Note that $\pp'$ has odd length. 
Let $A_1'=A_1\setminus V(\pp')$ and $B_1'=B_1\setminus V(\pp')$  
and we have $|A_1'| = |B_1'|=:m\in \frac k2 \mathbb N$. 
Further, by possibly extending $\pp$ by one $\e'$-good $(0,k/2)$-set and
 one $\e'$-good $(k/2,0)$-set we may assume that $\frac{2}{k}m$ is odd while conserving the above mentioned properties.
Moreover, we have that $|V(\pp')|\le |V(\pp)|+(b/4+2)\cdot k/2\le k^2\e_0' n+k$.

Next we pick $\e'$-good $(0, k/2)$-sets $L_1, L_2$, $\e'$-good $(k/2, 0)$-sets $S_1, S_2$ from $V\setminus V(\pp')$ such that $L' L_1$, $SS_1$, $L_2 S_2\in E(\h)$. 
Let $X$ be an arbitrary subset of $A_1'\setminus (S_1\cup S_2)$ of order $\frac{1}2(2m/k-1)$, and~$Y$ be an arbitrary subset of 
$B_1'\setminus (L_1\cup L_2)$ of order $\frac{1}2(2m/k-1)$.
We apply Lemma \ref{lem:finish2} on $\h[(B_1'\setminus Y)\cup X]$ with $\a_0=\sqrt\e_1$ and sets $L_1, L_2$ to obtain a Hamilton path $\pp_1$ with ends $L_1, L_2$. 
We apply Lemma \ref{lem:finish2} on $\h[(A_1'\setminus  X)\cup Y]$ with $\a_0=\sqrt\e_1$ and sets $S_1, S_2$ to obtain a Hamilton path $\pp_2$ with ends $S_1, S_2$. 
This yields a Hamilton $(k/2)$-cycle of $\h$
\[S~ \pp'~ L'~ L_1~ \pp_1~ L_2~ S_2~ \pp_2~ S_1~ S.\]

\medskip
\noindent\textit{Case $ii)$}. $k\in 4\mathbb{N}$.  

Then we have $b\in 2\mathbb N$ according to~\eqref{eq:bbbb}.
Note that $k/4$ and $k/4-1$ are two consecutive integers and without loss of generality, assume that $k/4$ is even. 
Now we greedily extend the path $\pp$ from $L$ by $b/2$ $\e'$-good $(k/4, k/4)$-sets and $b/2$ $\e'$-good $(k/4-1, k/4+1)$-sets. 
Note that the above process is possible since all edges involved are odd. 
Let the resulting path be $\pp'$ and note that it has odd length. Let $A_1'=A_1\setminus V(\pp')$ and $B_1'=B_1\setminus V(\pp')$  
and we have $|A_1'| = |B_1'|=:m\in \frac k2 \mathbb N$. 

Moreover, we have
\[
m - \frac{n-|V(\pp')|}{k}=|A_1\setminus V(\pp)| - \frac b2 \left(\frac k2 -1 \right) - \frac{n-|V(\pp)|-bk/2}{k}\equiv 0 - \frac b2 \cdot \frac k2 + b \equiv 0 \pmod 2,
\]
because $k/2, b\in 2\mathbb N$ and $|A_1\setminus V(\pp)|-\frac{n-|V(\pp)|}{k}$ is even by Claim \ref{clm:shortpath2}.
We claim that $m/k\in \mathbb N$. Indeed,
otherwise $2m/k$ must be odd, and from
\[
\frac{2m}{k} =\frac{|A_1'|+|B_1'|}{k} = \frac{n-|V(\pp')|}{k} \equiv m \pmod 2
\]
we infer that $m$ is odd.
Since $k/2$ is even, we obtain $2m/k\notin \mathbb N$, a contradiction.

Let the new end of $\pp'$ be $L'$ and
note that $|V(\pp')|\le |V(\pp)|+b\cdot k/2\le 4k^2\e_0' n$.
Next we pick disjoint $\e'$-good sets: $(0, k/2)$-sets $L_1, L_2$, $(k/2, 0)$-sets $S_1, S_2$, $(1,k/2-1)$-set $R_1$, and $(k/2-1,1)$-set $R_2$ from $V\setminus V(\pp')$ such that $
L' L_1, S R_1, R_1 S_1, L_2 R_2, R_2 S_2 \in E(\h).$ 
Let $X$ be an arbitrary subset of $A_1'\setminus (S_1\cup S_2\cup R_1\cup R_2)$ of order $m/k-1$, and $Y$ be an arbitrary subset of $B_1'\setminus (L_1\cup L_2\cup R_1\cup R_2)$ of order $m/k-1$.
Apply Lemma~\ref{lem:finish2} on $\h[(B_1'\setminus (Y\cup R_1\cup R_2))\cup X]$ with 
$\a_0=\sqrt\e_1$ and sets $L_1, L_2$ and obtain a Hamilton path $\pp_1$ with ends $L_1, L_2$. 
Then apply Lemma \ref{lem:finish2} again on $\h[(A_1'\setminus (X\cup R_1\cup R_2))\cup Y]$ 
with $\a_0=\sqrt\e_1$ and sets $S_1, S_2$ and obtain a Hamilton path $\pp_2$ with ends $S_1, S_2$. Thus, we get a Hamilton $(k/2)$-cycle
\[S~ \pp'~ L'~ L_1~ \pp_1~ L_2~ R_2~ S_2~ \pp_2~ S_1~ R_1~ S.\]

\subsubsection{The case $n\in \frac k2\mathbb N\setminus k \mathbb N$ and $\h$ is $\e$-close to ${\B}$}

We assume that the partition $(A_1, B_1)$ satisfies $\deg_{\B}(S)\le \deg_{\B}(T)$ for any $S\in \binom{A_1}{d}$ and any $T\in \binom{B_1}{d}$ -- if $(A_1, B_1)$ does not satisfy this, then swap $A_1$ and $B_1$.

Similar as in previous cases, we will prove the following claim.
\begin{claim}\label{clm:shortpath3}
There exists a $(k/2)$-path $\mathcal P$ in $\h$ such that 
\begin{itemize}
\item $|V(\pp)|\le 3k\e_0' n$,
\item $V(\pp)$ contains all $\e_1$-medium vertices,
\item the ends of $\mathcal P$ are  $\e'$-good $(k/2)$-sets with different parities,
\item $\pp$ has an even length and $|A_1\setminus V(\pp)|-\frac{n-|V(\pp)|}{k}$ is even.
\end{itemize}
\end{claim}

Before proving the claim we note that, the Step 2 is the same as  in the case $n\in k \mathbb N$ 
and $\h$ is $\e$-close to~${\B}$ since  $|A_1\setminus V(\pp)|-\frac{n-|V(\pp)|}{k}$ is even. 
It is thus left to prove Claim~\ref{clm:shortpath3}. 

\begin{proof}[Proof of Claim~\ref{clm:shortpath3}]
Let $\pp_{\mathrm b}$ be the path given by Lemma~\ref{clm:bridge3}. 
Let $\pp_M$ be the path with even ends given by Lemma~\ref{clm:pm} with $U=V(\pp_{\mathrm b})$, and then connect $\pp_M$ and $\pp_{\mathrm b}$ by one or two $(k/2)$-sets.
This is possible because their ends are $\e'$-good.
We extend the path by one more edge if its ends have the same parity and denote the resulting path by $\pp$. 
Note that $\pp$ contains one or three even edges and its two ends have different parities. 
This implies that $\pp$ has even length. We have $|V(\pp)|\le 3k\e_0' n$ similarly as in other cases. 
It remains to show that $|A_1\setminus V(\pp)|-\frac{n-|V(\pp)|}{k}$ is even.

\medskip
\noindent\textit{Case 1.} $\lfloor n/k \rfloor - |A_1|$ is odd. So $|A_1| \equiv n/k +1/2 \pmod 2$.
\medskip

In this case $\pp$ contains one even edge of the form 11 or three even edges of the form 00. In either case we have $|A_1\cap V(\pp)|\equiv \frac12 \left(\frac{|V(\pp)|}{k/2}+1 \right) \pmod 2$, namely, the number of digit 1's in $\pp$ and in a path of the form $11010\cdots 10$ are congruent modulo 2. Thus
\[
|A_1\setminus V(\pp)|-\frac{n-|V(\pp)|}{k}= |A_1| - \left( \frac nk +\frac12 \right)  - |A_1\cap V(\pp)|+  \left(\frac{|V(\pp)|}{k}+\frac12 \right),
\]
is even. 

\medskip
\noindent\textit{Case 2.} $\lfloor n/k \rfloor - |A_1|$ is even. So $|A_1| \equiv n/k -1/2 \pmod 2$.
\medskip

In this case $\pp$ contains one even edge of the form 00 or three even edges of the form 11. In either case we have $|A_1\cap V(\pp)|\equiv \frac12\left(\frac{|V(\pp)|}{k/2}-1\right) \pmod 2$, namely, the number of digit 1's in $\pp$ and in a path of the form $00101\cdots 01$ are congruent modulo 2. Thus
\[
|A_1\setminus V(\pp)|-\frac{n-|V(\pp)|}{k}= |A_1| - \left( \frac nk -\frac12 \right)  - |A_1\cap V(\pp)|+  \left(\frac{|V(\pp)|}{k}-\frac12 \right)
\]
is even. 
\end{proof}

\section{Deferred Proofs}\label{sec:deferred}

For a $k$-graph $\h$, let $\overline \h : = (V(\h), \binom{V(\h)}{k}\setminus E(\h))$.
To prove Lemmas~\ref{lem:finish2}--\ref{lem:finish1}, we need some results of Glebov, Person, and Weps \cite{GPW}. Given $1\le \ell \le k-1$ and $0\le \rho\le 1$, an ordered set $(z_1,\dots, z_{\ell})$ is 
$\rho$-\emph{typical} in a $k$-graph $\G$ if for every $i\in [\ell]$
\[
\deg_{\overline\G}( \{z_1, \dots, z_i\} ) \le \rho^{k-i}\binom {|V(\G)| - i}{k-i}.
\] 
It was shown in \cite{GPW} that every $k$-graph $\G$ with very large minimum vertex degree contains a tight Hamilton cycle. The proof of \cite[Theorem 2]{GPW} actually shows that 
any tight path of constant length with two typical ends can be extended to
a tight Hamilton cycle.  

\begin{theorem}[\cite{GPW}]
\label{thm:GPW}
Given $1\le \ell\le k-1$ and $0<\a \ll 1$, there exists an $m_0$ such that the following holds. 
Suppose that $\G$ is a $k$-graph on $V$ with $|V|= m\ge m_0$ and $\delta_1(\G)\ge (1 - \a)\binom{m-1}{k-1}$. 
Then given any two $(22\a)^{\frac1{k-1}}$-typical ordered $\ell$-sets $(z_1,\dots, z_{\ell})$ and $(y_1,\dots, y_{\ell})$, there exists a tight Hamilton path in $\G$ with ends
$(z_{\ell},z_{\ell-1},\dots, z_1)$ and  $(y_1,y_2,\dots, y_{\ell})$.\qed
\end{theorem}

We also use \cite[Lemma 3]{GPW}, in which $V^{2k-2}$ denotes the set of all $(2k-2)$-tuples  of not necessarily distinct elements of $V$.

\begin{lemma}[\cite{GPW}]
\label{lem:GPW}
Given $1\leq k$ and $0<\a \ll 1$, there exists an $m_0$ such that the following holds. 
Suppose that~$\G$ is a $k$-graph on $V$ with $|V|= m\ge m_0$ and $\delta_1(\G)\ge (1 - \a)\binom{m-1}{k-1}$. 
Then with probability at least ${8}/{11}$ a randomly selected $(x_1, \dots, x_{2k-2})\in V^{2k-2}$ 
is such that all $x_i$'s are  distinct and $(x_1, \dots, x_{k-1}), (x_{k}, \dots, x_{2k-2})$ are $(22\a)^{\frac1{k-1}}$-typical.\qed
\end{lemma}

Equipped with these auxiliary results we now give the proofs of Lemmas~\ref{lem:finish2}--\ref{lem:finish1}.

\begin{proof}[Proof of Lemma \ref{lem:finish}]

Let $\alpha_0\ll\alpha_1\ll\alpha\ll1$. 
By the assumption of the lemma we have $\delta_1(\h)\ge (1-\a_1)\binom{|Y|-1}{k-1}$. 
Let $Y'=Y\setminus (L_0\cup L_1)$ and $\h'=\h[Y']$. Since $t$ is large enough, we have,
\[
\delta_1(\h') \ge (1-\a_1)\binom{|Y|-1}{k-1} - k \binom{|Y|-1}{k-2} \ge (1-2\a_1)\binom{|Y'|-1}{k-1}.
\]
Since $L_0$ and $L_1$ are $\a_0$-good with respect to $\K^k(Y)$, we have $\deg_{\overline\h}(L_i)\le \a_0 |Y|^{k/2}\le \a_1 \binom{|Y|-k/2}{k/2}$ for $i\in\{0,1\}$. 
Thus, with probability at least $(1-\alpha)$, a random $k$-tuple $(z_1,\dots,z_{k/2},y_1,\dots, y_{k/1})\in (Y')^{k}$ 
satisfies
\begin{equation}\label{eq:L01}
L_0\cup\{z_1, \dots, z_{k/2}\}, \,L_1\cup\{y_1, \dots, y_{k/2}\}\in E(\h).
\end{equation}
Moreover, choosing 
 $(z_1, \dots, z_{k-1}, y_1, \dots, y_{k-1})$ from $(Y')^{2k-2}$ uniformly at random induces
 a uniform choice $(z_1,\dots,z_{k/2},y_1,\dots, y_{k/2})\in (Y')^{k}$. Thus, with $(1-\alpha)>3/11$,
 we infer from  Lemma~\ref{lem:GPW} that there exist $(44\a_1)^{\frac1{k-1}}$-typical  tuples 
 $(z_1, \dots, z_{k/2})$ and $(y_1, \dots, y_{k/2})$ of $k$ distinct elements for which~\eqref{eq:L01} holds.

By applying Theorem~\ref{thm:GPW} with $\ell=k/2$ and $\alpha$
we obtain a tight Hamilton path in $\h'$ with ends $(z_{k/2},z_{k/2-1},\dots, z_1)$  and  $(y_1,y_2,\dots, y_{k/2})$. 
Together with $L_0$ and $L_1$ this yields the desired $(k/2)$-path in $\h$.
\end{proof}

Next we prove Lemma \ref{lem:finish1}.
Another proof with a similar strategy can be found in \cite[Lemma 3.10]{HZ2}.

\begin{proof}[Proof of Lemma~\ref{lem:finish1}] Let $\alpha_0\ll\alpha\ll 1$ and let $X,Y, L_0,L_1$ be given. 
Throughout this proof  we refer to an $\alpha$-good set with respect to $\K_2^k(X,Y)$ simply as $\alpha$-good.
We call a set $S\subset V(\h)$ an $(a,b)$-set if $|S\cap X|=a$ and $|S\cap Y|=b$ and we further write $A B x$ for
$A\cup B\cup \{x\}$. 
Write $X=\{x_1, \dots, x_t\}$.
The main idea here is to partition $Y$ into $(k/2-1)$-sets $\{S_1,\dots, S_t\}$ where $S_1x_{r_1} S_2 x_{r_2}\cdots S_t x_{r_t}$ is the desired Hamilton path satisfying that $\{r_1, \dots, r_t\}=[t]$, $L_0=S_1\cup \{x_{r_1}\}$ and $L_1=S_t\cup \{x_{r_t}\}$.
To achieve this we plan to find a partition of $Y$ into $\{S_1,\dots, S_t\}$ so that we can carry out the following two-step process.
First for odd $i\in [t-1]$ we greedily choose $r_i$ such that $S_{i}x_{r_i}S_{i+1}$ (as a $(k-1)$-set) has high degree in $\h$ to the remaining vertices of $X$, denoted by $X'$.
Second we use Hall's Theorem on the auxiliary bipartite graph on $X'$ and $\{i\in [t-1], i $ even$\}$ where $x_j\sim i$ if and only if both $S_{i}x_{r_i}S_{i+1}x_j, S_{i+1} x_j S_{i+2}x_{r_{i+2}}\in E(\h)$.
Both of these would follow if all $S_i$'s are `typical', in particular, the second step requires that every vertex in $X'$ has high degree in $\Gamma$.
It is not clear to us how to find such a partition of $Y$, but it is easy to argue that the number of vertices in $X'$ with low degree is small.
To resolve this, before we partition $Y$ into $(k/2-1)$-sets, we first choose some random $(k/2-1)$-sets as buffer sets -- they can be used to match an arbitrary set of small number of vertices in $X$.
Namely, we will first set aside a small set $\A$ of structures in $X\cup Y$ and run the main idea on the remaining part of $\h$; then the small proportion of wrong vertices in $X'$ can be taken care by $\A$ so that we can apply Hall's Theorem on the majority of the bipartite graph $\Gamma$.

One main step in our proof is to establish the following claim.
\begin{claim}\label{claim:absorbpath2}
Let $q=4k\sqrt{\alpha}t$, then there  is a family $\A=\L\cup \R_{\mathrm{odd}}\cup \R_{\mathrm{even}}$  where $\L=\{L_1,L_2,\dots,L_{2q+1},L_1'\}$, 
$\R_{\mathrm{odd}}=\{R_{1}, R_3, R_5,\dots, R_{2q-1}\}$   are families of $(1,k/2-1)$-sets
and  $\R_{\mathrm{even}}=\{R_{2}, R_4,\dots, R_{2q}\}$ is a family of $(0,k/2-1)$-sets such that
\begin{itemize}
\item  $\A$ consists of  mutually disjoint subsets of $V\setminus L_0$ and $L_1'$ is $\alpha$-good. Furthermore,  $L_{2q+1}= L_1'$  if $|X|$ is odd and  $L_{2q+1}L_1'\in E(\h)$ if $|X|$ is even,
\item for each $R_{2i-1}\in\R_{\mathrm{odd}}$ we have $L_{2i-1}R_{2i-1}\in E(\h)$ and $R_{2i-1}L_{2i}\in E(\h)$,
\item for each $R_{2i}\in\R_{\mathrm{even}}$ the sets $L_{2i}R_{2i}$ and $R_{2i}L_{2i+1}$ are both $\alpha$-good and
conversely, for each $x\in X$ there are at least $q/2$ sets in $\R_{\mathrm{even}}$
such that $L_{2i}R_{2i}x\in E(\h)$ and $xR_{2i}L_{2i+1}\in E(\h)$.
\end{itemize}

\end{claim}
Note that once we have matched $\R_{\mathrm{even}}$  with a set $\{x_2,x_4,\dots,x_{2q}\}=X_{\mathrm{even}}\subset X\setminus V(\A)$  such that 
$L_{2i}R_{2i}x_{2i}, R_{2i}x_{2i}L_{2i+1}\in E(\h)$
we will obtain  the $k/2$-path 
\[L_1~ R_1~ L_2~ R_2~x_2~ L_3~ R_3~ L_4~R_4~x_4~L_5\cdots~ L_{2q-1}~ R_{2q-1}~ L_{2q}~ R_{2q}~x_{2q}~ L_{2q+1}~L_1'\]
with $\alpha$-good  ends $L_1$ and $L_1'$. Moreover, due to the third property of the claim we have the flexibility of choosing $X_{\mathrm{even}}$
to contain $q/2$ arbitrary vertices of $X$. Lastly, by the first property we can guarantee that $Y\setminus (V(\A)\cup L_0)$ is an odd multiple of $k/2$. 

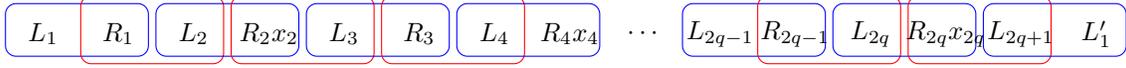
\begin{figure}[h]
\begin{center}
\begin{tikzpicture}
[inner sep=2pt,
   vertex/.style={circle, draw=blue!50, fill=blue!50},
   ]
\draw[color=blue, rounded corners] (-1,0) rectangle (0.9,0.7);
\draw[color=red, rounded corners] (0,-0.1) rectangle (1.9,0.8);
\draw[color=blue, rounded corners] (1,0) rectangle (2.9,0.7);
\draw[color=red, rounded corners] (2,-0.1) rectangle (3.9,0.8);
\draw[color=blue, rounded corners] (3,0) rectangle (4.9,0.7);
\draw[color=red, rounded corners] (4,-0.1) rectangle (5.9,0.8);
\draw[color=blue, rounded corners] (5,0) rectangle (6.9,0.7);
\node at (-0.5,0.3) {$L_1$};
\node at (0.5,0.3) {$R_1$};
\node at (1.5,0.3) {$L_2$};
\node at (2.5,0.3) {$R_2 x_2$};
\node at (3.5,0.3) {$L_3$};
\node at (4.5,0.3) {$R_3$};
\node at (5.5,0.3) {$L_4$};
\node at (6.5,0.3) {$R_4 x_4$};
\draw[color=blue, rounded corners] (8,0) rectangle (9.9,0.7);
\draw[color=red, rounded corners] (9,-0.1) rectangle (10.9,0.8);
\draw[color=blue, rounded corners] (10,0) rectangle (11.9,0.7);
\draw[color=red, rounded corners] (11,-0.1) rectangle (12.9,0.8);
\draw[color=blue, rounded corners] (12,0) rectangle (13.9,0.7);
\node at (8.5,0.3) {$L_{2q-1}$};
\node at (7.5,0.3) {$\cdots$};
\node at (9.5,0.3) {$R_{2q-1}$};
\node at (10.5,0.3) {$L_{2q}$};
\node at (11.5,0.3) {$R_{2q} x_{2q}$};
\node at (12.5,0.3) {$L_{2q+1}$};
\node at (13.5,0.3) {$L_1'$};
\end{tikzpicture}

\caption{The $(k/2)$-path formed by the members of $\A$}
\end{center}
\end{figure}

\begin{proof}[Proof of Claim~\ref{claim:absorbpath2}]For each $(2,\frac32k-3)$-set $S$ we fix a partition $S=LRL'$, where
$L, L'$ are $(1, k/2-1)$-sets  and $R$ is a $(0,\frac k2-1)$-set. 
Let $\mathcal F$ be the collection of  such sets $S=LRL'$, such that $L,L', LR$ and $RL'$ are all $\alpha$-good.
 Further, let $\S(x)$ denote the collection of those $S=LRL'$ such that $LRx,xRL'\in E(\h)$.
 We first establish that
 \begin{equation}\label{eq:FandS}
 |\F|\geq (1 - \a) \tbinom{|X|}2\tbinom{|Y|}{3k/2-3}\qquad\text{and}\qquad |\F\cap \S(x)|\geq (1 -2 \a) \tbinom{|X|}2\tbinom{|Y|}{3k/2-3}\text{ for all $x\in X$}.
 \end{equation}
To see \eqref{eq:FandS} let $B_{i,j}(x)$, with $i\leq 2$ and $j\leq k-2$ and $x\in X$,  denote the family of 
 $(i,j)$-sets in~$\h$, which contain $x$ and are
not ${\a}$-good.  
As all vertices are $\alpha_0$-good, 
 the number of $(2,k-1)$-sets in $\overline \h$ 
containing~$x$ is at most $\a_0 n^{k-1}$. On the other hand each element in $B_{i,j}(x)$ 
gives rise to at least $\alpha n^{k-i-j}$ such $(2,k-2)$-sets, thus 
\[
|B_{i,j}(x)|\leq \frac{\a_0 n^{k-1}}{\a n^{k-i-j}}\leq \a^2 n^{i+j-1}/2.
\]
By summing over $x$ we conclude that the number of $(1,\frac k2)$-sets in $\h$, which are not $\alpha$-good, is  at most $\alpha^2 n^{k/2}$ 
while the number of $(1,k-2)$-sets in $\h$, which are not $\alpha$-good, is at most $\alpha^2 n^{k-1}$.
Thus, the number of $(2,\frac 32k-3)$ sets $S=LRL'\not\in\F$ - i.e. such that some of the sets $L$, $L'$, $LR$, or $RL'$  is not $\alpha$-good - is at most
\[\alpha^2 n^{k/2}|X|\tbinom{|Y|}{k-2}+\alpha^2 n^{k-1}|X|\tbinom{|Y}{k/2-1}< \alpha\tbinom{|X|}{2}\tbinom{|Y|}{3k/2-3}.\]
This establishes the first part of~\eqref{eq:FandS}. Further, note that the
number of $(2,3k/2-3)$-sets $S\not\in \S(x)$ is at most $B_{1,k-2}(x)\cdot |X|\binom{|Y|}{k/2-1}\le \a \binom{|X|}2\binom{|Y|}{3k/2-3}$,
thus the second part follows from the first part.

We sequentially choose a family of $q$ elements from $\F$, making at each step a  random choice of an element disjoint from the chosen ones.
By \eqref{eq:FandS} and the fact that at most $k+q\cdot 3k/2$ vertices are already chosen, for any $x\in X$ and at each step, the probability that the random $(2, 3k/2-3)$-set is in $\S(x)$ is at least
\[
1-\frac{2\a \binom{|X|}2\binom{|Y|}{3k/2-3}}{\binom{|X^*|}2\binom{|Y^*|}{3k/2-3}} \ge 1-3\a \ge 3/4,
\]
where $X^*$ and $Y^*$ denote the intermediate sets of not chosen vertices in $X$ and $Y$, respectively.
By Lemma~\ref{lem:coupling} with $\delta = 1/3$ and the union bound, there exists a family $\F'=\{S_2, S_4, \dots, S_{2q}\}$ of disjoint $(2, \frac32k-3)$-sets, 
which contains at least $q/2$ members from each $\S(x)$, $x\in X$.
 Let $S_{2i}=L_{2i}R_{2i}L_{2i+1}$ be the implied partition of $S_{2i}$ which yields
the families  $\L=\{L_1,L_2,\dots,L_{2q+1}\}$ and
$\R_{\mathrm{even}}=\{R_{2}, R_4,\dots, R_{2q}\}$ with the required properties. We now choose a family $\R_{\mathrm{odd}}=\{R_{1},R_3,\dots,R_{2q-1}\}$
of disjoint $(0,\frac k2)$-sets from $V(\h)\setminus V(\F')$
such that  $L_{2i-1}R_{2i-1}, R_{2i-1}L_{2i}\in E(\h)$, $i\in[q]$, which is possible since $L_{2i-1}$ and $L_{2i}$ are both $\alpha$-good. 
Finally, if $t$ is even, we find an $\a$-good $(1,k/2-1)$-set $L_1'$ disjoint from $V(\A)$ such that $L_1'\cup L_{2q+1}\in E(\h)$. This is possible since $L_{2q+1}$ is $\a$-good. Otherwise, let $L_1':=L_{2q+1}$.
\end{proof}

Let $Y':=Y\setminus (V(\mathcal A)\cup L_0)$ and let $\G$ be the $(k-2)$-graph on $Y'$ which consists of all 
$\alpha$-good $(k-2)$-sets, i.e., which form an edge in $\h$ with all but at most 
$\alpha n^2$ elements from $\binom X2$. Then we have 
\[\delta_1(\G)\ge  (1-\a)\binom{|Y'|-1}{k-3},\]
since a $v\in Y'$ which violates this condition would be contained in at least $\alpha \binom{|Y'|-1}{k-3}\cdot \alpha n^2>\alpha n^{k-1}$ non-edges in $\h$, contradicting 
that $v$ is a $\alpha_0$-good vertex in $\h$.
Thus, by  Lemma \ref{lem:GPW} with probability at least $8/{11}$, uniformly chosen $(\hat z_1,\dots, \hat z_{k/2-1})$, $(\hat y_1,\dots, \hat y_{k/2-1})\in V^{k-1}$ form two disjoint ordered 
$(22\alpha)^{\frac1{k-1}}$-typical $(\frac k2-1)$-sets of $Y'$. 
Moreover the probability that $L_{1}'\cup \{\hat z_1,\dots, \hat z_{k/2-1}\}$ (or $L_0\cup \{\hat y_1,\dots, \hat y_{k/2-1}\}$, respectively) is  $\sqrt{\a}$-good is at least $1-\sqrt{\a}$ due to
the $\a$-goodness of $L_1'$ and $L_0$.
Therefore, there exists a choice  which satisfies both properties which we  denote by $(z_1,\dots, z_{k/2-1})$, $(y_1,\dots, y_{k/2-1})$.
Applying Theorem \ref{thm:GPW}  we obtain a tight Hamilton path of  $\G$ 
\[
\mathcal P=z_{k/2-1}z_{k/2-2}\cdots z_1 \cdots\cdots y_1y_2\cdots y_{k/2-1}
\]
and by following its order  we obtain a partition of $Y'$ into $(k/2-1)$-sets
\[
S_{1}=\{z_1,\dots, z_{k/2-1}\},\, S_{2},\dots, S_{t'}=\{y_1,\dots, y_{k/2-1}\}.
\]
Since $\mathcal P$ is a tight path in $\G$, we have that each $S_i S_{i+1}$ is $\a$-good $(0,k-2)$-set in $\h$ and by the choice of 
$(z_1,\dots, z_{k/2-1})$ and $(y_1,\dots, y_{k/2-1})$
from above we also have  that $L_{1}'\cup S_{1}, L_0\cup S_{t'}$ are $\sqrt{\a}$-good.
In the following we will match $\{S_i\}_{i\in[t']}\cup \R_{\mathrm{even}}$ 
with the vertices of $X\setminus V(\A)\cup L_0$ 
to form the desired Hamilton $(k/2)$-path of $\h$. To do so we will use the following two round process
\begin{enumerate}
\item Recall that $t'$ is odd.  We match $S_{2}, S_4,\dots, S_{t'-1}$  with suitable vertices $x_2,x_4,\dots, x_{t'-1}$ so that for each even $i\in[t']$
the sets $S_{i-1}S_{i} x_{i}$ and $S_{i} S_{i+1} x_{i}$ are $(\sqrt{\a}/2)$-good $(0,k-1)$-sets.
\item By making use of the properties of $\R_{\mathrm{even}}$ we then match the remaining vertices from $X$ with $\{S_1,S_3,\dots S_{t'}\}\cup \R_{\mathrm{even}}$ to complete
the Hamilton $(k/2)$-path in $\h$.\end{enumerate}

Concerning the first step we can simply greedily choose $x_2, x_4, \dots, x_{t'-1}\in X\setminus (V(\mathcal A)\cup L_0)$.
Note that   $|(V(\A)\cup L_1'\cup L_0)\cap X|\le 3q+2= 12k\sqrt\a t+2$. Moreover, as $S_{2i-1}S_{2i}$ and $S_{2i} S_{2i+1}$ are $\a$-good, 
there are at most $2\sqrt\a |V|\le k\sqrt{\a}t$ vertices $x\in X$ such that $S_{2i-1}S_{2i} x$ or $S_{2i} S_{2i+1} x$ is not $(\sqrt{\a}/2)$-good. 
Thus, all but at most $14k\sqrt\a t$ vertices in $X$ are not available as candidate for $x_{2i}$ initially and therefore the
process can be done greedily as we  only need $({t'-1})/2\le |X|/2$ vertices of $X$.

To carry out the second step  let $X_1:=X\setminus (V(\mathcal A)\cup L_1'\cup L_0\cup \{x_2, x_4, \dots, x_{t'-1}\})$ and note that $|X_1|= ({t'+1})/{2}+q$.
Let $I=\{1, 3, \dots, t'\}$
and consider the bipartite graph $\Gamma$ between $X_1$ and $I$ such that $x\in X_1$ is adjacent to an element $i\in I$ if and only if 
\begin{itemize}
\item $S_{i-1}x_{i-1}S_{i}x, S_{i}x S_{i+1}x_{i+1}\in E(\h)$ for odd $i\in [t']$,
\end{itemize}
where $S_0x_0=L_1'$ and $S_{t'+1}x_{t'+1}=L_0$.
Since both $S_{2i-1}S_{2i} x_{2i}$ and $S_{2i} S_{2i+1} x_{2i}$ are $(\sqrt{\a}/2)$-good, we have  $\deg_{\Gamma} (i) \ge |X_1|-2(\sqrt\a/2)n\ge |X_1|-k\sqrt \a t$ for every $i\in I$.
Let $X_0$ be the set of $x\in X_1$ such that $\deg_{\Gamma}(x)\le |I|/2$. Then
\[
|X_0| \frac{|I|}2 \le |X_1||I| - {e}_{\Gamma}(X_1, I) \le k \sqrt \a t\cdot |I|,
\]
which implies that $|X_0| \le 2k \sqrt\a t \leq q/2$. 
We match $\R_{\mathrm{even}}$ with a subset $X'_0\subset X_1$,  $X_0\subset X_0'$, matching to each 
$R_{2i}\in\R_{\mathrm{even}}$   an $x\in X_0'$ so that 
$L_{2i}R_{2i}x\in E(\h)$ and $xR_{2i}L_{2i+1}\in E(\h)$. This we do by first matching vertices from $X_0$ to elements of $\R_{\mathrm{even}}$  
and then matching the remaining members $\R_{\mathrm{even}}$ to elements in $X_1\setminus X_0$.
Due to the third property of Claim~\ref{claim:absorbpath2} this is possible.
Note that this completes $\A$ to a $(k/2)$-path with ends $L_1$ and $L_1'$ (see the remark after Claim~\ref{claim:absorbpath2}).

Finally let $X_2= X_1\setminus X_0'$ which has size $|X_2|=|X_1|-q=|I|$. Then $\Gamma'=\Gamma[X_2\cup I]$ is a graph with  $\delta(\Gamma')\ge |X_2|/2$.
Thus, by Hall's Theorem there is a perfect matching in $\Gamma'$ which gives the desired Hamilton path of $\h$.
\end{proof}

The proof of Lemma \ref{lem:finish2} is similar to the one from above, so we only give a sketch.

\begin{proof}[Proof sketch of Lemma \ref{lem:finish2}]
Let $\alpha_0\ll\alpha\ll 1$ and let $X,Y, L_0,L_1$ be given. 
Throughout this proof  we refer to an $\alpha$-good set with respect to $\K_1^k(X,Y)$ simply as $\alpha$-good.
Our goal is to write $X$ as $\{ x_1, \dots, x_{t} \}$ and partition $Y$ as
\[
\{L_1, R_1, L_2, R_2,\dots, L_t, R_t, L_{t+1} \}
\]
with $|L_i|=k/2$, $|R_i|=k/2-1$, and $L_0=L_{t+1}$ such that $L_i x_i R_i, x_i R_i L_{i+1} \in E(\h)$ for all $i\in [t]$.
Let $n:=|V(\h)|=kt+k/2$.

We first establish the following result.

\begin{claim}\label{claim:absorbpath3}
Let $q=12k\sqrt{\alpha}t$, then there  is a family $\A=\L\cup \R_{\mathrm{odd}}\cup \R_{\mathrm{even}}$  where $\L=\{L_1,L_2,\dots,L_{2q+1}\}$
is a family of $(0,k/2)$-sets,
$\R_{\mathrm{odd}}=\{R_{1}, R_3, R_5,\dots, R_{2q-1}\}$  a family of  $(1,k/2-1)$-sets
and  $\R_{\mathrm{even}}=\{R_{2}, R_4,\dots, R_{2q}\}$  a family of $(0,k/2-1)$-sets such that
\begin{itemize}
\item  $\A$ consists of  mutually disjoint subsets of $V\setminus L_0$ and $L_{2q+1}$ is $\alpha$-good. 
\item for each $R_{2i-1}\in\R_{\mathrm{odd}}$ we have $L_{2i-1}R_{2i-1}\in E(\h)$ and $R_{2i-1}L_{2i}\in E(\h)$,
\item for each $R_{2i}\in\R_{\mathrm{even}}$ the sets $L_{2i}R_{2i}$ and $R_{2i}L_{2i+1}$ are both $\alpha$-good and
conversely, for each $x\in X$ there are at least $q/2$ sets in $\R_{\mathrm{even}}$
such that $L_{2i}R_{2i}x\in E(\h)$ and $xR_{2i}L_{2i+1}\in E(\h)$.
\end{itemize}

\end{claim}

\begin{proof}[Proof sketch of Claim~\ref{claim:absorbpath3}]
For each $(0,\frac32k-1)$ set $S$ we fix a partition $S=LRL'$, where
$L, L'$ are $(0, k/2)$-sets  and $R$ is a $(0, k/2-1)$-set. 
Let $\mathcal F$ be the collection of  those sets $S=LRL'$, such that $L,L', LR$ and $RL'$ are all $\alpha$-good.
 Further, let $\S(x)$ denote the collection of those $S=LRL'$ such that $LRx,xRL'\in E(\h)$.
Similar to~\eqref{eq:FandS} we can establish that 
\begin{equation}\label{eq:FandS1}
|\F|\geq (1 - \a) \tbinom{|X|}2\tbinom{|Y|}{3k/2-3}\qquad\text{and}\qquad |\F\cap \S(x)|\geq (1 -2 \a) \tbinom{|X|}2\tbinom{|Y|}{3k/2-3}\text{ for all $x\in X$}.
\end{equation}
Then, by using Lemma~\ref{lem:coupling} we can pick a family 
$\F'=\{S_2, S_4,\dots, S_{2q}\}$ of pairwise disjoint sets from $\F$, 
which contains at least $q/2$ members from each $\S(x)$, $x\in X$.
 Let $S_{2i}=L_{2i}R_{2i}L_{2i+1}$ be the  partition of $S_{2i}$ which then yields
the families  $\L=\{L_1,L_2,\dots,L_{2q+1}\}$ and
$\R_{\mathrm{even}}=\{R_{2}, R_4,\dots, R_{2q}\}$ with the required properties. We now choose a family $\R_{\mathrm{odd}}=\{R_{1},R_3,\dots,R_{2q-1}\}$
of disjoint $(1, k/2)$-sets from $V(\h)\setminus V(\F')$
such that  $L_{2i-1}R_{2i-1}, R_{2i-1}L_{2i}\in E(\h)$, $i\in[q]$, which is possible since $L_{2i-1}$ and $L_{2i}$ are both $\alpha$-good. 
\end{proof}

Let $Y':=Y\setminus (V(\mathcal A)\cup L_0)$ and let $\G$ be the $(k-1)$-graph on $Y'$ which consists of all 
$\alpha$-good $(k-1)$-sets, i.e., which form an edge in $\h$ with all but at most 
$\alpha n$ elements from $X$. Then a similar calculation shows that 
\[\delta_1(\G)\ge  (1-\a)\binom{|Y'|-1}{k-3}\]
and following the approach in the previous proof, by Lemmas~\ref{lem:GPW} and~\ref{thm:GPW} we can find a partition of $Y_1$ as
\[
\{R_{2q+1},\, L_{2q+2}, \,R_{2q+2}, \dots, L_t, \, R_t\}.
\]
into  $(0,k/2-1)$-sets $R_i$'s and $(0,k/2)$-sets $L_i$'s such that both  $L_i R_i$ and $R_i L_{i+1}$ are $\sqrt\a$-good for $2q+1\le i\le t$, where $L_{t+1}=L_0$.

Let $X_1:=X\setminus V(\mathcal A)$ and note that $|X_1|= t-q$. 
Consider the bipartite graph $\Gamma$ between $X_1$ and $I=\{2q+1, 2q+2,\dots, t\}$ such that $x\in X_1$ is adjacent to an element $i\in I$ if and only if $L_i R_{i}x, R_i x L_{i+1}\in E(\h)$.
Since both of $L_i R_i$ and $R_i L_{i+1}$ are $\sqrt\a$-good, for every $i\in I$, $\deg_{\Gamma} (i) \ge |X_1| - 2\sqrt\a n\ge |X_1| - 3k\sqrt \a t$. 
Let $X_0$ be the set of $x\in X_1$ such that $\deg_{\Gamma}(x)\le |I|/2$. Then
\[
|X_0| \frac{|I|}2 \le |X_1||I| - {e}_{\Gamma}(X_1, I) \le 3k \sqrt \a t\cdot |I|,
\]
which implies that $|X_0| \le 6k \sqrt\a t $.
Thus we can match the vertices of $X$ to the structures in $Y$ and obtain a Hamilton path of $\h$ similar to the two steps in the proof of Lemma~\ref{lem:finish1}.
\end{proof}

\section{Concluding Remarks}

In this paper we found the minimum $d$-degree threshold for $(k/2)$-Hamiltonicity for all even $k\ge 6$ and all $d\ge k/2$. 
When $k=4$, we can add more edges to $\overline\B_{n,4}(A, B)$ and still avoid a Hamilton $2$-cycle. 
Partition $V$ into $A\cup B$ and fix two vertices $v_1, v_2\in A$. 
Let $\overline\B'_{n,4}(A, B)$ be the 4-graph obtained from $\overline\B_{n,4}(A, B)$ by adding all 4-sets $e$ of $V$ such that $|e\cap A|=3$ and $\{v_1, v_2\}\subseteq e$.

We claim that if $|A|$ is odd and $|A|\notin \{n/2, n/2+1\}$, then $\overline{\B}'_{n,4}(A,B)$ contains no Hamilton $2$-cycle.
Suppose to the contrary, that there is a Hamilton $2$-cycle $\C$ in $\overline{\B}'_{n,4}(A,B)$.
We represent $\C$ as a sequence of disjoint pairs of vertices $L_1, \dots, L_t$ with $t= n/2$. 
If all edges of $\C$ are even, then $L_i$'s are either all even or all odd.
Since $|A|=\sum_{i\in [t]}|A\cap L_i|$ is odd, all $L_i$'s must be odd, which implies that $|L_i\cap A|=1$. 
It follows that $|A|=n/2$, contradicting our assumption. 
Otherwise, $\C$ contains at least one odd edge.
However, since $\C$ is a cycle, $\C$ must contain an even number of odd edges (this can be seen by considering the parities of $L_i$'s).
By the definition of $\overline\B_{n,4}(A, B)$, $\C$ contains exactly two odd edges that both contain $\{v_1, v_2\}$.
We may thus assume that $L_1=\{v_1, v_2\}$ is even and all $L_i$, $i\neq 1$, are odd.
This implies that $|A|=n/2+1$, contradicting our assumption. 
Therefore we can add all $4$-graphs $\overline{\B}'_{n,4}(A,B)$ such that $|A|$ is odd and $|A|\notin \{n/2, n/2+1\}$ to $\h_{\ext}(n,4)$.

At last, we remark that for the missing case $k=4$ and $d=2$, by Theorem~\ref{lemNE}, it suffices to prove the extremal case, that is, find the best possible minimum $2$-degree condition for 2-Hamiltonicity in $4$-graphs which are close to either $\overline\B_{n,4}$ or $\B_{n,4}$.
For this case it is not clear to us how to build the `bridge', the short path overcoming the parity issue arising from the constructions (e.g., $\overline\B'_{n,4}(A, B)$) above, or whether our construction is indeed extremal.

\section*{Acknowledgement}
We thank two anonymous referees for their careful readings and helpful comments that improved the presentation of this paper.

\bibliographystyle{abbrv}
\bibliography{Bibref}

\begin{thebibliography}{10}

\bibitem{ABHKP16}
P.~Allen, J.~B\"ottcher, H.~H{\`a}n, Y.~Kohayakawa, and Y.~Person.
\newblock Blow-up lemmas for sparse graphs.
\newblock {\em ArXiv e-prints}, 2016.

\bibitem{BMSSS1}
J.~O. Bastos, G.~O. Mota, M.~Schacht, J.~Schnitzer, and F.~Schulenburg.
\newblock Loose {H}amiltonian cycles forced by large $(k-2)$-degree -
  approximation version.
\newblock {\em SIAM J. Discrete Math.}, 31:2328--2347, 2017.

\bibitem{BMSSS2}
J.~O. Bastos, G.~O. Mota, M.~Schacht, J.~Schnitzer, and F.~Schulenburg.
\newblock Loose {H}amiltonian cycles forced by large $(k-2)$-degree - sharp
  version.
\newblock {\em Contributions to Discrete Mathematics}, 13(2):88--100, 2019.

\bibitem{BHS}
E.~Bu{\ss}, H.~H{\`a}n, and M.~Schacht.
\newblock Minimum vertex degree conditions for loose {H}amilton cycles in
  3-uniform hypergraphs.
\newblock {\em J. Combin. Theory Ser. B}, 103(6):658--678, 2013.

\bibitem{CzMo}
A.~Czygrinow and T.~Molla.
\newblock Tight codegree condition for the existence of loose {H}amilton cycles
  in 3-graphs.
\newblock {\em SIAM J. Discrete Math.}, 28(1):67--76, 2014.

\bibitem{EKR}
P.~Erd{\H{o}}s, C.~Ko, and R.~Rado.
\newblock Intersection theorems for systems of finite sets.
\newblock {\em Quart. J. Math. Oxford Ser. (2)}, 12:313--320, 1961.

\bibitem{FrFu85}
P.~Frankl and Z.~F\"{u}redi.
\newblock Forbidding just one intersection.
\newblock {\em Journal of Combinatorial Theory, Series A}, 39(2):160 -- 176,
  1985.

\bibitem{GaMy}
F.~Garbe and R.~Mycroft.
\newblock Hamilton cycles in hypergraphs below the dirac threshold.
\newblock {\em Journal of Combinatorial Theory, Series B}, 133:153--210, 2018.

\bibitem{GPW}
R.~Glebov, Y.~Person, and W.~Weps.
\newblock On extremal hypergraphs for {H}amiltonian cycles.
\newblock {\em European J. Combin.}, 33(4):544--555, 2012.

\bibitem{HS}
H.~H\`an and M.~Schacht.
\newblock Dirac-type results for loose {Hamilton} cycles in uniform
  hypergraphs.
\newblock {\em J. Combin. Theory Ser. B}, 100:332--346, 2010.

\bibitem{HZ2}
J.~Han and Y.~Zhao.
\newblock Minimum codegree threshold for hamilton $\ell$-cycles in $k$-uniform
  hypergraphs.
\newblock {\em J. Combin. Theory Ser. A}, 132(0):194 -- 223, 2015.

\bibitem{HZ1}
J.~Han and Y.~Zhao.
\newblock Minimum degree thresholds for loose {Hamilton} cycle in 3-graphs.
\newblock {\em J. Combin. Theory Ser. B}, 114:70 -- 96, 2015.

\bibitem{HZ_forbidHC}
J.~Han and Y.~Zhao.
\newblock Forbidding {H}amilton cycles in uniform hypergraphs.
\newblock {\em J. Combin. Theory Ser. A}, 143:107 -- 115, 2016.

\bibitem{HM67}
A.~J.~W. Hilton and E.~C. Milner.
\newblock Some intersection theorems for systems of finite sets.
\newblock {\em Quart. J. Math. Oxford Ser. (2)}, 18:369--384, 1967.

\bibitem{KK}
G.~Katona and H.~Kierstead.
\newblock Hamiltonian chains in hypergraphs.
\newblock {\em J. Graph Theory}, 30(2):205--212, 1999.

\bibitem{KKMO}
P.~Keevash, D.~K\"uhn, R.~Mycroft, and D.~Osthus.
\newblock Loose {Hamilton} cycles in hypergraphs.
\newblock {\em Discrete Math.}, 311(7):544--559, 2011.

\bibitem{KMO}
D.~K\"uhn, R.~Mycroft, and D.~Osthus.
\newblock Hamilton $\ell$-cycles in uniform hypergraphs.
\newblock {\em J. Combin. Theory Ser. A}, 117(7):910--927, 2010.

\bibitem{KuOs14ICM}
D.~K{\"u}hn and D.~Osthus.
\newblock Hamilton cycles in graphs and hypergraphs: an extremal perspective.
\newblock {\em Proceedings of the International Congress of Mathematicians
  2014, Seoul, Korea}, Vol 4:381--406, 2014.

\bibitem{MaRu}
K.~Markstr\"{o}m and A.~Ruci\'{n}ski.
\newblock Perfect {M}atchings (and {H}amilton {C}ycles) in {H}ypergraphs with
  {L}arge {D}egrees.
\newblock {\em European J. Comb.}, 32(5):677--687, July 2011.

\bibitem{RRRSS}
C.~Reiher, V.~R{\"o}dl, A.~Ruci{\'n}ski, M.~Schacht, and E.~Szemer{\'e}di.
\newblock Minimum vertex degree condition for tight hamiltonian cycles in
  3-uniform hypergraphs.
\newblock {\em Proc. London Math. Soc.}, 119:409--439, 2019.

\bibitem{RR}
V.~R\"odl and A.~Ruci\'nski.
\newblock Dirac-type questions for hypergraphs --- a survey (or more problems
  for endre to solve).
\newblock {\em An Irregular Mind}, Bolyai Soc. Math. Studies 21:561--590, 2010.

\bibitem{RRS06}
V.~R\"odl, A.~Ruci\'nski, and E.~Szemer\'edi.
\newblock A {D}irac-type theorem for 3-uniform hypergraphs.
\newblock {\em Combin. Probab. Comput.}, 15(1-2):229--251, 2006.

\bibitem{RRS08}
V.~R\"odl, A.~Ruci\'nski, and E.~Szemer\'edi.
\newblock An approximate {D}irac-type theorem for $k$-uniform hypergraphs.
\newblock {\em Combinatorica}, 28(2):229--260, 2008.

\bibitem{RRS09}
V.~R{\"o}dl, A.~Ruci{\'n}ski, and E.~Szemer{\'e}di.
\newblock Perfect matchings in large uniform hypergraphs with large minimum
  collective degree.
\newblock {\em J. Combin. Theory Ser. A}, 116(3):613--636, 2009.

\bibitem{RRS11}
V.~R\"odl, A.~Ruci\'nski, and E.~Szemer\'edi.
\newblock Dirac-type conditions for {Hamiltonian} paths and cycles in 3-uniform
  hypergraphs.
\newblock {\em Advances in Mathematics}, 227(3):1225--1299, 2011.

\bibitem{Sze}
E.~Szemer{\'e}di.
\newblock Regular partitions of graphs.
\newblock In {\em Probl\`emes combinatoires et th\'eorie des graphes ({C}olloq.
  {I}nternat. {CNRS}, {U}niv. {O}rsay, {O}rsay, 1976)}, volume 260 of {\em
  Colloq. Internat. CNRS}, pages 399--401. CNRS, Paris, 1978.

\bibitem{TrZh12}
A.~Treglown and Y.~Zhao.
\newblock Exact minimum degree thresholds for perfect matchings in uniform
  hypergraphs.
\newblock {\em J. Combin. Theory Ser. A}, 119(7):1500--1522, 2012.

\bibitem{TrZh13}
A.~Treglown and Y.~Zhao.
\newblock Exact minimum degree thresholds for perfect matchings in uniform
  hypergraphs {II}.
\newblock {\em J. Combin. Theory Ser. A}, 120(7):1463--1482, 2013.

\bibitem{zsurvey}
Y.~Zhao.
\newblock Recent advances on dirac-type problems for hypergraphs.
\newblock In {\em Recent Trends in Combinatorics}, volume 159 of {\em the IMA
  Volumes in Mathematics and its Applications}. Springer, New York, 2016.

\end{thebibliography}

\end{document}